\newcommand{\R}{I\!\!R}
\newcommand{\Z}{\mathbb Z}
\newcommand{\N}{I\!\!N}
\newcommand{\dist}{\mathrm{dist}}

\newcommand{\Dim}{\mathrm{dim}}

\newcommand{\second}{\mathrm I\!\mathrm I}

\newcommand{\Ddt}{\tfrac{\mathrm D}{\mathrm dt}}
\newcommand{\Dds}{\tfrac{\mathrm D}{\mathrm ds}}

\newcommand{\Ddtt}{\tfrac{\mathrm D^2}{\mathrm dt^2}}

\documentclass[twoside,letterpaper,10pt]{amsart}

\usepackage{epsfig}
\usepackage{times}
\usepackage{hyperref}

\numberwithin{equation}{section}

\title[Homoclinic Orbits in Riemannian Manifolds]%
{Orthogonal Geodesic Chords, Brake Orbits and Homoclinic Orbits in
Riemannian Manifolds}

\author[R. Giamb\`o]{Roberto Giamb\`o}
\author[F. Giannoni]{Fabio Giannoni}
\author[P. Piccione]{Paolo Piccione}
\address{Dipartimento di Matematica e Informatica,\hfill\break\indent
Universit\`a di Camerino, Italy} \email{roberto.giambo@unicam.it,
fabio.giannoni@unicam.it, \hfill\break\phantom{\sl  E-mail
address: } \ \ paolo.piccione@unicam.it}

\urladdr{http://www.ime.usp.br/\~{}piccione}

\thanks{The third author's permanent address is: Departamento de Matem\'atica, Instituto de
Matem\'atica e Estat\'\i stica, Universidade de S\~ao Paulo,
Brazil.}


\date{March 4th, 2004}

\begin{document}


\theoremstyle{plain}\newtheorem{teo}{Theorem}[section]
\theoremstyle{plain}\newtheorem{prop}[teo]{Proposition}
\theoremstyle{plain}\newtheorem{lem}[teo]{Lemma}
\theoremstyle{plain}\newtheorem{cor}[teo]{Corollary}
\theoremstyle{definition}\newtheorem{defin}[teo]{Definition}
\theoremstyle{remark}\newtheorem{rem}[teo]{Remark}
\theoremstyle{definition}\newtheorem{example}[teo]{Example}

\theoremstyle{plain}\newtheorem*{teon}{Theorem}


\begin{abstract}
 The study of solutions with fixed energy
 of certain classes of Lagrangian
 (or Hamiltonian) systems is reduced, via the classical Maupertuis--Jacobi
 variational principle, to the study of geodesics in Riemannian
 manifolds.  We are interested in investigating the
 problem of existence of brake orbits and homoclinic orbits,
 in which case the Maupertuis--Jacobi principle produces a Riemannian
 manifold with boundary and with metric degenerating in a non trivial way
 on the boundary.   In this paper we use  the
 classical Maupertuis--Jacobi principle to show how to remove the
 degeneration of the metric on the boundary,
 and we prove in full generality how the brake orbit and the homoclinic
 orbit multiplicity problem can be reduced to the study  of multiplicity
 of orthogonal geodesic chords in a manifold with {\em regular\/} and {\em strongly concave\/}
 boundary.
\end{abstract}

\maketitle

\tableofcontents

\begin{section}{Introduction}
\label{sec:intro}

The study of periodic and homoclinic orbits of Lagrangian and
Hamiltonian systems is an extremely  active research field
in classical and modern mathematics, having a huge number of
applications in physical sciences. One of the peculiarities
of the problem is that, although already very popular among
classical analysts and geometers, it has never been out of
fashion, and it has been studied along the time with techniques
of an increasing level of sophistication. Indeed, the study of
solutions of Hamiltonian systems has motivated many recent
developments of several mathematical theories, including
Calculus of Variations, Symplectic Geometry and Morse Theory,
among others, and the vaste literature on the topic witnesses
the leading role of the subject in modern mathematics.

The central interest of the present paper is to study solutions of
an autonomous Lagrangian (or Hamiltonian) system, having
prescribed energy, in a manifold $M$ that belong to two special
classes of solutions: the homoclinic orbits and the brake orbits.
Homoclinic orbits are solutions $x:\R\to M$ of the system for
which the limits $\lim\limits_{t\to+\infty}x(t)$ and
$\lim\limits_{t\to-\infty}x(t)$  exist and are equal, and
$\lim\limits_{t\to\pm\infty}\dot x(t)=0$. Such limits must then be
a critical point of the potential function of the system.
Brake orbits are a special
class of periodic solutions that have an oscillating character,
i.e., periodic solutions $x:\R\to M$ having period $2T$, with
$x(T+t)=x(T-t)$ and $\dot x(T+t)=-\dot x(T-t)$ for all $t\in\R$.
Clearly, $\dot x(kT)=0$ for all $k\in\Z$.

By a classical variational principle, known  as the
Maupertuis--Jacobi principle, solutions of autonomous Lagrangian
or Hamiltonian systems having a fixed value of the energy
correspond to geodesics relatively to a Riemannian metric, called
the Jacobi metric. When dealing with  homoclinic orbits issuing
from a critical point of the potential function, or with brake
orbits, then the classical formulation of the Maupertuis--Jacobi
principle fails, due to the fact that such solutions pass through
a region where the Jacobi metric degenerates in a non trivial way.
An accurate analysis of the geodesic behavior near such
degeneracies, that occur on the boundary of the level set of the
potential function, has lead many authors to obtain existence
results by perturbation techniques. More specifically, following
an original idea by Seifert \cite{seifert}, some authors (see \cite{gluckziller})
have been able to
perform a geometrical construction consisting in attaching a
smooth, {\em convex\/} and sufficiently small collar (see
Figure~\ref{fig:torocollo})
\begin{figure}
\begin{center}
\psfull \epsfig{file=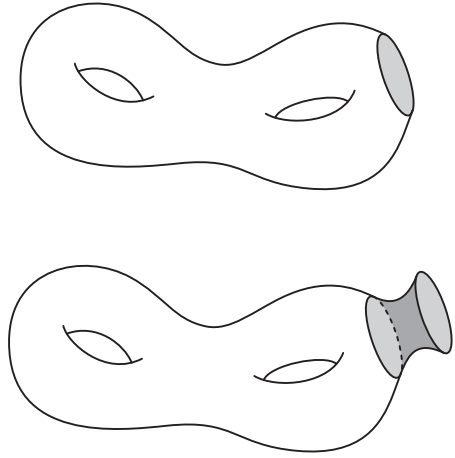, height=6cm} \caption{Gluing
a collar with convex boundary to a concave
boundary.}\label{fig:torocollo}
\end{center}
\end{figure}to the
degenerate region, in such a way that the geodesics in the
resulting manifold could be counted by standard techniques in
convex Riemannian geometry (\cite{bos, LustSchn}). Then, a limit argument was used to
obtain existence results for geodesics in the original degenerate
metric by letting the size of the collar go to zero.  The same
idea cannot be used if one wants to obtain multiplicity results,
due to the fact that such limit procedure does not guarantee that
possibly distinct geodesics in the perturbed metric converge to
geometrically distinct geodesics in the original Jacobi metric,
unless one poses {\em ad hoc\/} "non resonance" assumptions (see \cite{gluckziller}).
Here, by geometrically distinct, we mean geodesics having different
images; the non resonance assumptions mentioned above guarantees
that it is avoided the situation in which distinct geodesics in the
perturbed metric tend to the same periodic geodesic travelled
a different number of times.

The starting point of this paper is the idea that, if one wants
to preserve the number of distinct geodesics, then one has to
perform a geometrical construction
that avoids limits procedure. Such construction would obviously be
based on a careful investigation of the geodesic
behavior near the boundary of the level set of the potential
function. Working in this direction has lead to the quite
remarkable observation that the boundary of a non critical
level set of the potential function, or of a small ball around
a non degenerate maximum point of the potential, are near
certain hypersurfaces that are {\em strongly concave\/} relatively
to the Jacobi metric, and that have the property that orthogonal
geodesic chords arriving on one of these hypersurfaces can be uniquely
extended to geodesic chords up to the degenerate boundary.
The presence of concave hypersurfaces near the degenerate boundary
can be interpreted as an indication that Seifert's
technique of gluing a convex collar would be somewhat innatural in order
to study the multiplicity problem in full generality.

The main results of this paper are contained in Theorem~\ref{thm:lem8.5},
relating the brake orbits problem to the orthogonal geodesic chords
problem, and Theorem~\ref{thm:propfin}, that deals with the homoclinics problem.

 The issue of concavity, as opposed to the {\em convexity\/}
 property used in the classical literature, is the key point to
 develop a multiplicity theory for brake orbits and homoclinic
 orbits under purely topological assumptions on the underlying
 manifolds. These multiplicity results constitute the topic
 of two forthcoming papers by the authors (\cite{GiaGiaPicI, GiaGiaPicII}).
\end{section}

\begin{section}{Geodesics and Concavity}
\label{sec:geodes} Let $(M,g)$ be a smooth (i.e., of class $C^2$)
Riemannian manifold with $\Dim(M)=m\ge2$, let $\dist$ denote the
distance function on $M$ induced by $g$; the symbol $\nabla$ will
denote the covariant derivative of the Levi-Civita connection of
$g$, as well as the gradient differential operator for smooth maps
on $M$. The Hessian $\mathrm H^f(q)$ of a smooth map $f:M\to\R$ at
a point $q\in M$ is the symmetric bilinear form $\mathrm
H^f(q)(v,w)=g\big((\nabla_v\nabla f)(q),w\big)$ for all $v,w\in
T_xM$; equivalently, $\mathrm H^f(q)(v,v)=\frac{\mathrm
d^2}{\mathrm ds^2}\big\vert_{s=0} f(\gamma(s))$, where
$\gamma:\left]-\varepsilon,\varepsilon\right[\to M$ is the unique
(affinely parameterized) geodesic in $M$ with $\gamma(0)=q$ and
$\dot\gamma(0)=v$. We will denote by $\Ddt$ the covariant
derivative along a curve, in such a way that $\Ddt\dot x=0$ is the
equation of the geodesics. A basic reference on the background material
for Riemannian geometry is \cite{docarmo}.

Let $\Omega\subset M$ be an open subset;
$\overline\Omega=\Omega\bigcup\partial \Omega$ will denote its
closure. There are several notion of convexity and concavity in
Riemannian geometry, extending the usual ones for subsets of the
Euclidean space $\R^m$. In this paper we will use a somewhat 
concavity assumption for compact subsets of $M$, that we
will refer as "strong concavity" below,
and which is stable by
$C^2$-small perturbations of the boundary. Let us first recall the
following:
\begin{defin}\label{thm:defconvexity}
$\overline\Omega$ is said to be {\em convex\/} if every geodesic
 $\gamma:[a,b]\to \overline\Omega$ whose endpoints
$\gamma(a)$ and $\gamma(b)$ are in $\Omega$ has image entirely
contained in $\Omega$. Likewise, $\overline\Omega$ is said to be
{\em concave\/} if its complement $M\setminus\overline\Omega$ is
convex.
\end{defin}

If $\partial \Omega$ is a smooth embedded submanifold of $M$, let
$\second_{\mathfrak n}(x):T_x(\partial\Omega)\times
T_x(\partial\Omega)\to\R$ denote the {\em second fundamental form
of $\partial\Omega$ in the normal direction $\mathfrak n\in
T_x(\partial\Omega)^\perp$}. Recall that $\second_{\mathfrak
n}(x)$ is a symmetric bilinear form on $T_x(\partial\Omega)$
defined by:
\[\phantom{\qquad v,w\in T_x(\partial\Omega),}\second_{\mathfrak n}(x)(v,w)=g(\nabla_vW,
\mathfrak n),\qquad v,w\in T_x(\partial\Omega),\] where $W$ is any
local extension of $w$ to a smooth vector field along
$\partial\Omega$.

\begin{rem}\label{thm:remphisecond}
Assume that it is given a smooth function $\phi:M\to\R$ with the
property that $\Omega=\phi^{-1}\big(\left]-\infty,0\right[\big)$
and $\partial\Omega=\phi^{-1}(0)$, with $\mathrm d\phi\ne0$ on
$\partial\Omega$. \footnote{For example one can choose $\phi$ such that
$\vert\phi(q)\vert=\dist(q,\partial\Omega)$ for all $q$ in a
(closed) neighborhood of $\partial\Omega$.} The following equality
between the Hessian $\mathrm H^\phi$ and the second fundamental
form\footnote{%
Observe that, with our definition of $\phi$, then $\nabla\phi$ is
a normal vector to $\partial\Omega$ pointing {\em outwards\/} from
$\Omega$.} of $\partial\Omega$ holds:
\begin{equation}\label{eq:seches}
\phantom{\quad x\in\partial\Omega,\ v\in
T_x(\partial\Omega);}\mathrm H^\phi(x)(v,v)=
-\second_{\nabla\phi(x)}(x)(v,v),\quad x\in\partial\Omega,\ v\in
T_x(\partial\Omega);\end{equation} Namely, if
$x\in\partial\Omega$, $v\in T_x(\partial\Omega)$ and $V$ is a
local extension around $x$ of $v$ to a vector field which is
tangent to $\partial\Omega$, then $v\big(g(\nabla\phi,V)\big)=0$
on $\partial\Omega$, and thus:
\[\mathrm H^\phi(x)(v,v)=v\big(g(\nabla\phi,V)\big)-g(\nabla\phi,\nabla_vV)=-\second_{\nabla\phi(x)}(x)(v,v).\]

Note that the second fundamental form is defined intrinsically,
while there is general no natural choice for a function $\phi$
describing the boundary of $\Omega$ as above.
\end{rem}

\begin{defin}\label{thm:defstrongconcavity}
We will say that that $\overline\Omega$
is {\em strongly concave\/} if $\second_{\mathfrak n}(x)$ is
positive definite for all $x\in\partial\Omega$ and all inward
pointing normal direction $\mathfrak n$.
\end{defin}

\begin{rem}\label{thm:remopencondition}
Strong concavity is evidently a {\em $C^2$-open
condition}. It should also be emphasized that if $\overline\Omega$
is strongly concave, then for {\em any\/} smooth map $\phi:M\to\R$
as in Remark~\ref{thm:remphisecond}, then for all $q\in\partial\Omega$,
the Hessian $\mathrm H^\phi(q)$ is
negative definite on $T_q\big(\partial\Omega\big)$.
From this observation, it follows immediately that geodesics starting
tangentially to $\partial\Omega$ move inside $\Omega$.
\end{rem}

The main objects of our study are geodesics in $M$ having image in
$\overline\Omega$ and with endpoints orthogonal to
$\partial\Omega$. We distinguish a special class of such geodesics,
called "weak", whose relevance will not be emphasized in the present paper,
but it will be used in a substantial way in the proof of the multiplicity results in
\cite{GiaGiaPicI, GiaGiaPicII}.
\begin{defin}\label{thm:defOGC}
A geodesic $\gamma:[a,b]\to M$ is called a {\em geodesic chord\/}
in $\overline\Omega$ if
$\gamma\big(\left]a,b\right[\big)\subset\Omega$ and
$\gamma(a),\gamma(b)\in\partial\Omega$; by a {\em weak geodesic
chord\/} we will mean a geodesic $\gamma:[a,b]\to M$ with image in
$\overline\Omega$ and endpoints
$\gamma(a),\gamma(b)\in\partial\Omega$. A (weak) geodesic chord is
called {\em orthogonal\/} if $\dot\gamma(a^+)\in
(T_{\gamma(a)}\partial\Omega)^\perp$ and $\dot\gamma(b^-)\in
(T_{\gamma(b)}\partial\Omega)^\perp$, where
$\dot\gamma(\,\cdot\,^\pm)$ denote the lateral derivatives (see
Figure~\ref{fig:WOGC}). 
An orthogonal geodesic chord in $\overline\Omega$ whose
endpoints belong to distinct connected components of $\partial\Omega$
will be called a {\em crossing orthogonal geodesic chord\/} in $\overline\Omega$.
\end{defin}
\begin{figure}
\begin{center}
\psfull \epsfig{file=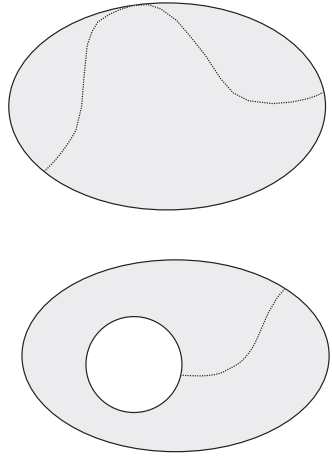, height=6cm} \caption{A weak
orthogonal geodesic chord (WOGC) in $\overline\Omega$
(above), and a crossing OGC (below).}\label{fig:WOGC}
\end{center}
\end{figure}

For shortness, we will write \textbf{OGC} for ``orthogonal
geodesic chord'' and \textbf{WOGC} for ``weak orthogonal geodesic
chord''.

For the proof of the multiplicity results in \cite{GiaGiaPicI, GiaGiaPicII},
we will use a geometrical construction that will work in a situation where one can
exclude {\em a priori\/} the existence in $\overline\Omega$ of (crossing)
weak orthogonal geodesic chords in $\partial\Omega$.
We will now show that one does not lose generality in assuming
that there are no such WOGC's in $\overline\Omega$ by proving the
following:
\begin{prop}
\label{thm:noWOGC}
Let $\Omega\subset M$ be an open set
whose boundary $\partial\Omega$ is smooth and compact and with $\overline\Omega$
strongly concave.
Assume that there are only a finite number of  (crossing) orthogonal
geodesic chords in $\overline\Omega$. Then, there exists an
open subset $\Omega'\subset\Omega$ with the following properties:
\begin{enumerate}
\item\label{itm:nowogcs1} $\overline{\Omega'}$ is diffeomorphic to $\overline\Omega$
and it has smooth boundary;
\item\label{itm:nowogcs2} $\overline{\Omega'}$ is strongly concave;
\item\label{itm:nowogcs3} the number of (crossing) OGC's in $\overline{\Omega'}$  is less than
or equal to the number of (crossing) OGC's in $\overline\Omega$ ;
\item\label{itm:nowogcs4} every (crossing) WOGC in $\overline{\Omega'}$
is a (crossing) OGC in $\overline{\Omega'}$.
\end{enumerate}
\end{prop}
\begin{proof}
The desired set $\Omega'$ will be taken of the form:
\[\Omega'=\phi^{-1}\big(\left]-\infty,-\delta\right[\big),\]
with $\delta>0$ small, and with $\phi$ a smooth map as in Remark~\ref{thm:remphisecond}
such that $\vert\phi(q)\vert=\dist(q,\partial\Omega)$ for $q$ near $\partial\Omega$.
Observe that if $\delta$ is small enough,
then by continuity $\mathrm d\phi\ne0$ on $\phi^{-1}([-\delta,0])$,
which implies that $\partial\Omega'$ is smooth and that $\overline{\Omega'}$ is diffeomorphic
to $\overline\Omega$, as we see using the integral curves of $\nabla\phi$.
Since strong concavity is an open condition
in the $C^2$-topology, if $\delta>0$ is small enough then
$\overline{\Omega'}$ is strongly concave, proving \eqref{itm:nowogcs2}.

Moreover, $\delta$ must be chosen small enough so that
the exponential map gives a diffeomorphism from
an open neighborhood of the zero section of the normal bundle of
$\partial\Omega$ to the set $\phi^{-1}\big(\left]-2\delta,2\delta\right[\big)$;
the existence of such $\delta$ is guaranteed by our compactness
assumption on $\partial\Omega$.
Since $\phi(q)=-\dist(q,\partial\Omega)$ near $\partial\Omega$, then
every (crossing) geodesic in $\Omega'$ that arrives orthogonally at
$\partial\Omega'$ can be smoothly extended to a (crossing) geodesic in $\overline\Omega$
that arrives orthogonally at $\partial\Omega$; observe that
any such extended geodesic only touches $\partial\Omega$ at the
endpoints, i.e., it is a (crossing) OGC in $\overline\Omega$. This  proves
part \eqref{itm:nowogcs3}.

We claim that there exists $\delta>0$ arbitrarily small such that
every (crossing) WOGC  is a (crossing) OGC in $\phi^{-1}\big(\left]-\infty,-\delta\right]\big)$.
Assume on the contrary that there exists a sequence $\delta_n>0$
with $\delta_n\to0$ as $n\to\infty$, a sequence $0 <s_n< 1$ and
a sequence of (crossing) geodesics $\gamma_n:[0,1]\to\Omega$ with
$\phi(\gamma_n(0))=\phi(\gamma_n(s_n))=\phi(\gamma_n(1))=-\delta_n$,
$\dot\gamma_n(0)$ and $\dot\gamma_n(1)$ orthogonal to $\phi^{-1}(-\delta_n)$
and $\phi(\gamma_n(s))\le-\delta_n$ for all $s\in[0,1]$ and all $n\in\N$.
As we have observed, for $n$ large each geodesic $\gamma_n$ can be smoothly
extended to a (crossing) OGC in $\overline\Omega$, and clearly
all such extensions cannot make a {\em finite\/} set of
geometrically distinct (crossing) OGC's in $\overline\Omega$.
Namely, each $\gamma_n$ is tangent to the surface
$\phi^{-1}(-\delta_n)$, and to {\em no other\/} surface
of the form $\phi^{-1}(-\delta)$ with $\delta<\delta_n$.
This says that the extensions of the $\gamma_n$ are
all geometrically distinct, which
contradicts the fact that there is only a finite number of (crossing)
OGC's in $\overline\Omega$ and proves part \eqref{itm:nowogcs4}.
\end{proof}

\end{section}
\begin{section}{Brake and Homoclinic Orbits of Hamiltonian Systems}
\label{sec:brakehom} Let $p=(p_i)$, $q=(q^i)$ be coordinates on
$\R^{2m}$, and let us consider a {\em natural\/} Hamiltonian
function $H\in C^2\big(\R^{2m},\R\big)$, i.e., a function of the
form
\begin{equation}\label{eq:hamfun}
H(p,q)=\frac 12 \sum_{i,j=1}^m a^{ij}(q)p_ip_j+V(q),
\end{equation}
where $V\in C^2\big(\R^{m},\R\big)$ and $A(q)=\big(a^{ij}(q)\big)$
is a positive definite quadratic form on $\R^m$:
\[\sum_{i,j=1}^m a^{ij}(q)p_ip_j\ge\nu(q)\vert q\vert^2\]
for some continuous function $\nu:\R^m\to\R^+$ and for all
$(p,q)\in\R^{2m}$.

The corresponding Hamiltonian system is:
\begin{equation}\label{eq:HS}
\left\{
\begin{aligned}
&\dot p=-\frac{\partial H}{\partial q}\\
&\dot q=\frac{\partial H}{\partial p},
\end{aligned}
\right.
\end{equation}
where the dot denotes differentiation with respect to time.

For all $q\in\R^m$, denote by $\mathcal L(q):\R^m\to\R^m$ the
linear isomorphism whose matrix with respect to the canonical
basis is $\big(a_{ij}(q)\big)$, the inverse of
$\big(a^{ij}(q)\big)$; it is easily seen that, if $(p,q)$ is a
solution of class $C^1$ of \eqref{eq:HS}, then $q$ is actually  a
map of class $C^2$ and
\begin{equation}\label{eq:pintermsofq}
p=\mathcal L(q)\dot q.
\end{equation} With a
slight abuse of language, we will say that a $C^2$-map $q:I\to
\R^m$ is a solution of \eqref{eq:HS} if $(p,q)$ is a solution of
\eqref{eq:HS} where $p$ is given by \eqref{eq:pintermsofq}. Since
the system \eqref{eq:HS} is autonomous, i.e., time independent,
then the function $H$ is constant along each solution, and it
represents the total energy of the solution of the dynamical
system. There exists a large amount of literature concerning the
study of periodic solutions of autonomous Hamiltonian systems
having energy $H$ prescribed (see for instance \cite{Long} and the
references therein).

We will be concerned with a special kind of periodic solutions of
\eqref{eq:HS}, called {\em brake orbits}. A brake orbit for the
system \eqref{eq:HS} is a non constant periodic solution $\R\ni
t\mapsto\big(p(t),q(t)\big) \in\R^{2m}$ of class $C^2$ with the
property that $p(0)=p(T)=0$ for some $T>0$. Since $H$ is
even in the variable $p$, a brake orbit $(p,q)$ is $2T$-periodic,
with $p$ odd and $q$ even about $t=0$ and about $t=T$. Clearly, if
$E$ is the energy of a brake orbit $(p,q)$, then
$V\big(q(0)\big)=V\big(q(T)\big)=E$.

The link between solutions of brake orbits and orthogonal geodesic
chords is obtained in Theorem~\ref{thm:lem8.5} (used in \cite{GiaGiaPicII} to obtain
the multiplicity result for brake orbits). Its proof is based on
a well known variational principle, that
relates solutions of \eqref{eq:HS} having prescribed energy $E$
with curves in the open subset $\Omega_E\subset\R^m$:
\begin{equation}\label{eq:jacobispace}
\Omega_E=V^{-1}\big(\left]-\infty,E\right[\big)=\big\{x\in\R^m:V(x)<E\big\}\end{equation}
endowed with the {\em Jacobi metric} (see Proposition
\ref{thm:MauJac}):
\begin{equation}\label{eq:jacobimetric}
g_E(x)=\big(E-V(x)\big)\cdot\frac 12\sum_{i,j=1}^m
a_{ij}(x)\,\mathrm dx^i\,\mathrm dx^j.
\end{equation}
Let us now consider the problem of homoclinics on a Riemannian manifold
$(M,g)$.

Assume that we are given a map $V\in C^2\big(M,\R\big)$; the
corresponding second order Hamiltonian system is the equation:
\begin{equation}\label{eq:secondorder}
\Ddt\dot q+\nabla V(q)=0.
\end{equation}
Note that if $M=\R^m$ and $g$ is the Riemannian metric
\begin{equation}\label{eq:Riemmetric}
g=\frac12\sum_{i,j=1}^m a_{ij}(x)\,\mathrm dx^i\,\mathrm dx^j,
\end{equation}
where the coefficients $a_{ij}$ are as above, then equation
\eqref{eq:secondorder} is equivalent to \eqref{eq:HS}, in the
sense that $x$ is a solution of \eqref{eq:secondorder} if and only
if the pair $q=x$ and $p=\mathcal L(x)\dot x$ is a solution of
\eqref{eq:HS}.

Let $x_0\in M$ be a critical point of $V$, i.e., such that $\nabla
V(x_0)=0$. We recall that a {\em homoclinic orbit\/} for the system
\eqref{eq:secondorder} emanating from $x_0$ is a solution $q\in
C^2\big(\R,M\big)$ of \eqref{eq:secondorder} such that:
\begin{eqnarray}
&\label{eq:1.5}
\lim\limits_{t\to-\infty}q(t)=\lim\limits_{t\to+\infty}q(t)=x_0,
\\
&\label{eq:1.6} \lim\limits_{t\to-\infty}\dot
q(t)=\lim\limits_{t\to+\infty}\dot q(t)=0.
\end{eqnarray}
To the authors' knowledge, the only result available in the literature
on multiplicity of homoclinics in the autonomous case is due to
Ambrosetti and Coti--Zelati \cite{AmbCot}, to Rabinowitz \cite{Rab} and to Tanaka \cite{Tan}.
A quite general multiplicity result for homoclinics, generalizing
those in \cite{AmbCot} and in \cite{Tan}, will be given in \cite{GiaGiaPicI} using
the result of Theorem~\ref{thm:propfin}.

It should also be mentioned that very likely all the results
in this paper can be extended to the case of  Hamiltonian
functions $H$ more general than \eqref{eq:hamfun}. As observed
by Weinstein in \cite{Wein}, Hamiltonians that are positively homogeneous in
the momenta lead to Finsler metrics rather than Riemannian metrics.

\end{section}
\begin{section}{The Maupertuis Principle}
\label{sec:notations}
Throughout this section, $(M,g)$ will denote a Riemannian manifold of
class $C^2$; all our constructions will be made in suitable
(relatively) compact subsets of $M$, and for this reason it will
not be restrictive to assume, as we will, that $(M,g)$ is
complete.

\subsection{The variational framework}
\label{sub:varframe}

The symbol $H^1\big([a,b],\R^m\big)$ will denote the Sobolev space
of all absolutely continuous function $f:[a,b]\to\R^m$ whose weak
derivative is square integrable. Similarly, $H^1\big([a,b],M\big)$
will denote the infinite dimensional Hilbert manifold consisting
of all absolutely continuous curves $x:[a,b]\to M$ such that
$\varphi\circ x\vert_{[c,d]}\in H^1\big([c,d],\R^m)$ for all chart
$\varphi:U\subset M\to \R^m$ of $M$ such that $x\big([c,d]\big)
\subset U$. By $H^1_{\textrm{loc}}\big(\left]a,b\right[,\R^m\big)$
we will denote the vector space of all continuous maps
$f:\left]a,b\right[\to\R^m$ such that $f\vert_{[c,d]}\in
H^1\big([c,d],\R^m\big)$ for all $[c,d]\subset\left]a,b\right[$;
the set $H^1_{\textrm{loc}}\big(\left]a,b\right[,M\big)$ is
defined similarly. The Hilbert space norm of
$H^1\big([a,b],\R^m\big)$ will be denoted by
$\Vert\cdot\Vert_{a,b}$; for the purposes of this paper it will
not be necessary to make the choice among equivalent norms of
$H^1\big([a,b],\R^m\big)$.

\subsection{The Maupertuis--Jacobi principle for brake orbits}\label{sub:Mau-Jac}
Let $V\in C^2\big(M,\R\big)$ and let $E\in\R$.
Consider  the sublevel $\Omega_E$ of $V$ in
\eqref{eq:jacobispace} and the {\em Maupertuis integral\/}
$f_{a,b}:H^1\big([a,b],\Omega_E\big)\to\R$, which is the geodesic
action functional relative to the metric $g_E$
\eqref{eq:jacobimetric}, given by:
\begin{equation}\label{eq:maupint}
f_{a,b}(x)=\frac12\int_a^b\big(E-V(x)\big)g\big(\dot x,\dot
x\big)\, \mathrm dt,
\end{equation}
where $g$ is the Riemannian metric \eqref{eq:Riemmetric}. Observe
that the metric $g_E$ degenerates on $\partial\Omega_E$.

The functional $f_{a,b}$ is smooth, and its differential is
readily computed as:
\begin{equation}\label{eq:diffmaupint}
\mathrm df_{a,b}(x)W=\int_a^b\big(E-V(x)\big)g\big(\dot x,\Ddt
W\big)\, \mathrm dt-\frac12\int_a^bg\big(\dot x,\dot x\big)
g(\nabla V(x),W\big)\,\mathrm dt,
\end{equation}
where $W\in H^1\big([a,b],\R^m\big)$. The corresponding
Euler--Lagrange equation of the critical points of $f_{a,b}$ is
\begin{equation}\label{eq:pr1}
\big(E-V(x(s))\big)\Ddt\dot x(s) - g\big(\nabla V(x(s)),\dot
x(s)\big)\dot x(s) +\frac 12 g\big(\dot x(s),\dot x(s)\big)\nabla
V(x(s))=0,
\end{equation}
for all $s\in ]a,b[$.

Solutions of the Hamiltonian system \eqref{eq:HS} having fixed
energy $E$ and critical points of the functional $f_{a,b}$ of
\eqref{eq:maupint} are related by the following variational
principle, known in the literature as the {\em Maupertuis--Jacobi
principle}:
\begin{prop}
\label{thm:MauJac} Assume that $E$ is a regular value of the
function $V$.

Let $x\in C^0\big([a,b],\R^m\big)\cap
H^1_{\textrm{loc}}\big(\left]a,b\right[,\R^m\big)$ be a non
constant curve such that
\begin{equation}\label{eq:i}
\int_a^b\big(E-V(x)\big)g\big(\dot x,\Ddt W\big)\, \mathrm
dt-\frac12\int_a^bg\big(\dot x,\dot x\big) g(\nabla
V(x),W\big)\,\mathrm dt=0
\end{equation}
for all $W\in C^\infty_0\big(\left]a,b\right[,\R^m\big)$, and such
that:
\begin{equation}\label{thm:eq:ii}
V\big(x(s)\big)<E,\quad\text{for all $s\in\left]a,b\right[$;}
\end{equation}
and
\begin{equation}\label{eq:iii}
V\big(x(a)\big),V\big(x(b)\big)\le E.
\end{equation}
Then, $x\in H^1\big([a,b],\R^m\big)$, and if
$\,\,V\big(x(a)\big)=V\big(x(b)\big)=E$, it is $x(a)\ne x(b)$.
Moreover, in the above situation, there exist positive constants
$c_x$ and $T$ and a $C^1$-diffeo\-morphism $\sigma:[0,T]\to[a,b]$
such that:
\begin{equation}\label{eq:3p2}
\big(E-V(x)\big)g\big(\dot x,\dot x\big)\equiv c_x\quad \text{on
$[a,b]$},
\end{equation}
and, setting $q=x\circ\sigma:[0,T]\to\R^m$, and $p(s)=\mathcal
L(q(s))\dot q(s)$, the pair $(q,p):[0,T]\to\R^{2m}$ is a solution
of \eqref{eq:HS} having energy $E$ with $q(0)=x(a)$, $q(T)=x(b)$.
If $V\big(x(a)\big)=V\big(x(b)\big)=E$ then $q$ can be extended to
a $2T$-periodic brake orbit of \eqref{eq:HS}.
\end{prop}
\begin{proof}
A proof when $\mathcal L$ is the identity map $\mathrm{id}$ can be
found for instance in \cite{benci}. For convenience of the reader
we give here a sketch of the proof in the general case.

Since $x$ satisfies \eqref{eq:i}, standard regularization
arguments show that $x$ is of class $C^2$ on $]a,b[$, while
integration by parts gives \eqref{eq:pr1} $\forall s\in]a,b[$.
Equation \eqref{eq:3p2} follows contracting both sides of
\eqref{eq:pr1} with $\dot x$ using $g$. Now set
\begin{equation}\label{eq:pr2}
t(s)=\frac12\int_a^s\frac{c_x}{E-V(x(\tau))}\,\text d\tau.
\end{equation}
A simple estimate shows that $T\equiv t(b)<+\infty$. Indeed,
setting $$C=\sup\{g\big(\nabla V(x),\nabla
V(x)\big)^{1/2}\,:\,x\in\overline\Omega_E\},$$ and using
\eqref{eq:3p2}, one has
\[
\left|\frac{\text d}{\text
ds}\left(\frac{1}{E-V(x(s))}\right)\right|\le \frac{C\, g\big(\dot
x,\dot x\big)^{1/2}}{\big(E-V(x)\big)^2}=\frac{C\sqrt{c_x}}
{\big(E-V(x)\big)^{5/2}}.
\]
Therefore, standard estimates for ordinary differential equations
gives the existence of a constant $D_x$ such that
\[
\frac{1}{E-V(x(s))}\le
D_x\left(\frac{1}{(s-a)^{2/3}}+\frac{1}{(b-s)^{2/3}}\right),\qquad\forall
s\in ]a,b[,
\]
proving that $t(b)<+\infty$ and that $x\in H^1([a,b],\R^m)$.

Now, denote by $\sigma:[0,T]\to [a,b]$ the inverse map of
\eqref{eq:pr2}, and set $q(t)=x(\sigma(t))$. Since $\sigma'(t)=2
(c_x)^{-1} \big(E-V(x(\sigma(t)))\big)$, a straightforward
computation shows that $\Dds\dot q=-\nabla V(q)$ and $\frac 12
g\big(\dot q,\dot q\big)+V(q)\equiv E$. Therefore, the pair
$(q,\mathcal L(q)\dot q):[0,T]\to\R^{2m}$ is a solution of
\eqref{eq:HS} with energy $E$.

Moreover $q(0)=x(a)$ and $q(T)=x(b)$, and by the uniqueness of the
Cauchy problem, if $V(x(a))=V(x(b))=E$ it must be $q(0)\not=q(T)$,
and $q$ can be extended to a periodic brake orbit.
\end{proof}

\subsection{The Maupertuis--Jacobi Principle near a nondegenerate
maximum of the potential energy.} \label{sub:Jacsing}
The above formulation of the Maupertuis--Jacobi principle is not
suited to study homoclinic orbits issuing from a critical point
of the potential function $V$. Our next goal is to establish an
extension of the principle that will be applied in this situation.

\begin{prop}\label{thm:MauJacnearmax}
Let $(M,g)$ be a Riemannian manifold, $V\in C^2\big(M,\R\big)$, let
$x_0\in M$ be a nondegenerate maximum of $V$, and set $E=V(x_0)$. Assume that
$x$ is a curve in the set $ C^0\big([a,b],\overline
\Omega_E\big)\bigcap
H^1_{\textrm{loc}}\big(\left[a,b\right[,\overline\Omega_E\big)$
such that:
\begin{equation}\label{eq:2.4bis}
\int_a^b\big(E-V(x)\big)g\big(\dot x,\Ddt W\big)\, \mathrm
dt-\frac12\int_a^bg\big(\dot x,\dot x\big) g(\nabla
V(x),W\big)\,\mathrm dt=0
\end{equation}
for all $W\in C^\infty_0\big(\left]a,b\right[,\R^m\big)$, and such
that
\begin{eqnarray}
\label{eq:2.5bis} &&V\big(x(s)\big)<E,\ \text{for $s\in\left[a,b\right[$};\\
\label{eq:2.6bis} && x(b)=x_0.
\end{eqnarray}
Then, there exists a $C^1$-diffeomorphism
$\sigma:\left[0,+\infty\right[\to\left[a,b\right[$ such that the
curve $q=x\circ\sigma$ is a solution of \eqref{eq:secondorder}
satisfying $q(0)=x(a)$ and $\lim\limits_{t\to+\infty}q(t)=x_0$,
$\lim\limits_{t\to+\infty}\dot q(t)=0$.
\end{prop}
\begin{proof}
Choose $\varrho\in\left]0,\dist\big(x(a),x_0\big)\right[$ and
define $\alpha_1\in\left]a,b\right[$ as the {\em first\/} instant
$s$ at which $\dist\big(x(s),x_0\big)=\varrho$. By
\eqref{eq:2.4bis}, the restriction $x\vert_{[a,\alpha_1]}$ is a
geodesic relatively to the metric $g_E$, since
$x\big([a,\alpha_1]\big)$ is contained in a region where $E-V$ is
positive. Denote by $c_x$ the constant value of $(E-V(x))g(\dot
x,\dot x)$; for all $s\in[a,\alpha_1]$ set:
\[
t(s)=\frac 12\int_a^s\frac{c_x}{E-V(x(\tau))}\,\mathrm d\tau
\]
and denote by $\sigma:[0,t(\alpha_1)]\to[a,\alpha_1]$ the inverse
function of $s\mapsto t(s)$. Then, a straightforward calculations
shows that the map $q=x\circ\sigma$ is a solution of the equation
\eqref{eq:secondorder} with $\frac12g(\dot q,\dot q)+V(q)\equiv E$
on $[0,s(\alpha_1)]$.

Let us choose $\alpha_2\in\left]\alpha_1,b\right[$ be the {\em
first\/} instant $s$ at which $\dist(x(s),x_0)=\frac\varrho2$; we
can repeat the construction above obtaining  a solution $q_*$ of
\eqref{eq:secondorder} defined on an interval $[0,t(\alpha_2)]$.
The key observation here is that, in fact, such a function $q_*$
is an extension of $q$,  and therefore it satisfies the same
conservation law $\frac12g(\dot{q_*}, \dot {q_*})+V(q_*)\equiv E$
on $[0,t(\alpha_2)]$. An iteration of this construction produces a
sequence $a<\alpha_1<\alpha_2<\ldots<b$ such that
$\dist(x(\alpha_k),x_0)=\frac\varrho{2^{k-1}}$, maps of class
$C^1$, $t:\left[a,L\right[\to\left[0,T\right[$, its inverse
$\sigma:\left[0,T\right[\to \left[a,L\right[$, where:
\[
T=\frac 12\int_a^{L}\!\!\!\frac{c_x}{E-V(x(\tau))}\,\mathrm
d\tau\in\left]0,+\infty\right], \quad
L=\lim_{k\to\infty}\alpha_k\in\left]a,b\right],
\]
and a curve of class $C^2$, $q=x\circ\sigma:\left[0,T\right[\to
\overline\Omega_E$, that satisfies \eqref{eq:secondorder}, and
with
\begin{equation}\label{eq:conslawdim}\frac12g(\dot{q},
\dot {q})+V(q)\equiv E
\end{equation}   on $\left[0,T\right[$; in particular,
$g(\dot q,\dot q)$ is bounded.

Let us prove that $T=+\infty$ and that $\lim\limits_{t\to
+\infty}q(t)=x_0$. We know that, by construction,
$\lim\limits_{k\to\infty}t(\alpha_k)=T$ and
$\lim\limits_{k\to\infty}q(t(\alpha_k))=x_0$; suppose by absurd
that there exists $\bar\rho>0$, and a sequence $\beta_k$ such that
$\lim\limits_{k\to\infty}\beta_k=L$ and
$\dist(q(t(\beta_k)),x_0)\ge\bar\rho$ for all $k$. Since $x_0$ is
an isolated maximum point, we can assume $\bar\rho$ small enough
so that
\begin{equation}\label{eq:maxV}
\inf\limits_{\frac12\bar\rho\le\dist(Q,x_0)\le\bar\rho}
\big(E-V(Q)\big)\equiv\bar e>0.
\end{equation}
Up to subsequences, we can obviously assume that
$\beta_k\in\left]\alpha_k, \alpha_{k+1}\right]$ for all $k$; for
$k$ sufficiently large, there exists
$\gamma_k\in\left]\alpha_k,\beta_k\right[$ which is the first
instant $t\in \left]\alpha_k,\beta_k\right[$ at which
$\dist\big(q(s(t)),x_0\big)=\frac{\bar\rho}2$. Since $g(\dot
q,\dot q)$ is bounded, there exists $\bar\nu>0$ such that
\begin{equation}\label{eq:abs1}
t(\gamma_k)-t(\alpha_k)\ge\bar\nu,\quad \text{for all $k$;}
\end{equation}
from \eqref{eq:maxV} and \eqref{eq:abs1} we get:
\begin{equation}\label{eq:intdivergente}
\int_0^{t(\alpha_{N+1})}\!\!\!\!\!\!\!\! \!\!\!\!
\big(E-V(q(\tau))\big)\,\mathrm d\tau\ge
\sum_{k=1}^N\int_{t(\alpha_k)}^{t(\gamma_k)}\!\!\!\!\!
\big(E-V(q(\tau))\big)\,\mathrm d\tau\ge \sum_{k=1}^N \bar e
\bar\nu=N\bar e\bar\nu \longrightarrow+\infty
\end{equation}
 as $N\to\infty$. On the other hand,
for all $s\in\left]a,L\right[$,
\[\int_0^{t(s)}\big(E-V(q(\tau))\big)\,\mathrm d\tau=
\frac 12\int_a^s c_x\;\mathrm d\theta= \frac{(b-a)}2 c_x,\] which
is obviously inconsistent with \eqref{eq:intdivergente}, and
therefore proves that $\lim\limits_{t\to T^-}q(t)=x_0$. Moreover,
the conservation law \eqref{eq:conslawdim} implies that
$\lim\limits_{t\to T^-}\dot q(t)=0$.

Finally, the local  uniqueness of the solution of an initial value
problem implies immediately that $T$ cannot be finite; for, the
only solution $q$ of \eqref{eq:secondorder} satisfying $q(T)=x_0$
and $\dot q(T)=0$ is the constant $q\equiv x_0$.
\end{proof}

\end{section}
\begin{section}%
{Orthogonal Geodesic Chords and the Maupertuis Integral.} \label{sec:mapertuis} In this section we will
prove the main result of the paper, showing how to reduce the brake
orbit and the homoclinics multiplicity problem to a multiplicity
result for orthogonal geodesic chords.

We will begin with the study of the Jacobi metric near the level surface
$V^{-1}(E)$, with $E$ regular value of $V$.

\subsection{The Jacobi distance near a regular value of the potential.}
Let $g$ be a Riemannian metric, $g_E=\big(E-V(x)\big) g$,
$\Omega_E$ as in \eqref{eq:jacobispace}; assume $\nabla
V(x)\not=0$ for all $x\in V^{-1}(E)$ and that $\overline{\Omega}_E$ is
compact.

\begin{lem}\label{thm:lem8.2}
For all $Q\in\Omega_E$, the infimum:
\[d_E(Q)\!:=\!\inf\Big\{\!\!\int_0^1\!\!\!\!\!\left((E-V(x))g\big(\dot x,\dot x\big)
\right)^{1/2}\!\!\mathrm dt:x\!\in\!
H^1\big([0,1],\overline\Omega_E\big),\ x(0)\!=\!Q,\ x(1)
\in\partial\Omega\Big\}\] is attained on at least one curve
$\gamma_Q\in H^1\big([0,1],\overline\Omega_E\big)$ such that
$\big(E-V(\gamma_Q)\big)g\big(\dot\gamma_Q,\dot\gamma_Q\big)$ is
constant, $\gamma_Q\big(\left[0,1\right[\big)\subset\Omega$, and
$\gamma_Q$ is a $C^2$ curve on $\left[0,1\right[$. Moreover, such a
curve satisfies assumption \eqref{eq:i} of
Proposition~\ref{thm:MauJac} on the interval $[a,b]=[0,1]$.
\end{lem}
\begin{proof}
For all $k\in\N$ sufficiently large, set
$\Omega_k=V^{-1}\big(\left]-\infty,E-\frac1k\right[\big)\subset\Omega_E$,
and consider the problem of minimization of the $g_E$-length
functional:
\[L_E(x)=\int_0^1\big[(E-V(x))g(\dot x,\dot x)\big]^{\frac12}\,\mathrm ds,\]
in the space $\mathfrak G_k$ consisting of curves $x\in
H^1\big([0,1],\overline\Omega_k\big)$ with $x(0)=Q$ and
$x(1)\in\partial\Omega_k$.

It is not hard to prove, by standard arguments, that for all
$\Omega_k\not=\emptyset$, the above problem has a solution
$\gamma_k$ which is a $g_E$-geodesic, and with
$\gamma_k\big(\left[0,1\right[\big)\subset \Omega_k$.

Set $q_k=\gamma_k(1)\in\partial\Omega_k$ and $l_k=L_E(\gamma_k)$.
Since $q_k$ approaches $\partial\Omega$ as $k\to\infty$, arguing
by contradiction we get:
\[\liminf_{k\to\infty} l_k\ge d_E(Q).\]
Now, if by absurd it was:
\[\liminf_{k\to\infty} l_k> d_E(Q),\]
then we could find a curve $x\in
H^1\big([0,1],\overline\Omega\big)$ with $x(0)=Q$,
$x(1)\in\partial\Omega$, and with
$L_E(x)<\liminf\limits_{k\to\infty} l_k$. 
Then, a suitable
reparameterization of $x$ would yield a curve $y\in\mathfrak G_k$
with $L_E(y)<l_k$, which contradicts the minimality of $l_k$ and
proves that
\begin{equation}\label{eq:eq10.11}
\liminf_{k\to\infty}l_k=d_E(Q).
\end{equation}
Now, arguing as in the proof of Proposition~\ref{thm:MauJac}, we
see that the sequence:
\begin{equation}\label{eq:seq}
\int_0^1\frac{\mathrm dt}{E-V(\gamma_k(t))}
\end{equation}
is bounded. Now, $\int_0^1\big
(E-V(\gamma_k)\big)g\big(\dot\gamma_k,\dot\gamma_k\big)\,\text
d\tau=l_k^2\equiv(E-V(\gamma_k)\big)g\big(\dot\gamma_k,\dot\gamma_k\big)$
is bounded, which implies $\int_0^1
g\big(\dot\gamma_k,\dot\gamma_k\big)\,\text d\tau$ bounded, namely
the sequence $\gamma_k$ is bounded in
$H^1\big([0,1],\overline\Omega_E\big)$. Up to subsequences, we
have a curve $\gamma_Q\in H^1\big([0,1],\overline\Omega_E\big)$
which is an $H^1$-weak limit of the $\gamma_k$'s; in particular,
$\gamma_k$ is uniformly convergent to $\gamma_Q$.

We claim that such a curve $\gamma_Q$ satisfies the required
properties. First, $\gamma_Q([0,1[)\subset\Omega_E$. Otherwise, if
$b<1$ is the first instant where $\gamma_Q(b)\in\partial\Omega_E$,
by \eqref{eq:eq10.11} and the conservation law of the energy for
$\gamma_k$ one should have
\[
(b-1)l_k^2=\int_b^1\big(E-V(\gamma_k)\big)g\big(\dot\gamma_k,\dot\gamma_k\big)\,\text
d\tau\longrightarrow 0,
\]
in contradiction with $Q\not\in\partial\Omega_E$. Then $\gamma_Q$
satisfies \eqref{eq:i} in $[0,1]$ since it is a $H^1$--weak limit
of $\gamma_k$, which is a sequence of $g_E$--geodesics.

Clearly, $\gamma_Q$ is of class $C^2$ on $\left[0,1\right[$,
because the convergence on each interval $\left[0,b\right]$ is
indeed smooth for all $b<1$.

Finally, since $L_E(z)\le\liminf\limits_{k\to\infty}l_k$, from
\eqref{eq:eq10.11} it follows that $L_E(\gamma_Q)=d_E(Q)$, and this
concludes the proof.
\end{proof}

\begin{rem}\label{rem:rem3.3}
It is immediate to see that, $\gamma_Q$ is a minimizer as in Lemma
\ref{thm:lem8.2} if and only if is a minimizer for the functional
\begin{equation}\label{eq:funct01}
f_{0,1}(x)=\frac 12\int_0^1\big(E-V(x)\big)g\big(\dot x,\dot
x\big)\,\text dt
\end{equation}
in the space of curves
\begin{equation}\label{eq:Xspace}
X_Q=\{x\in H^1([0,1],\overline\Omega_E)\,:\,x(0)=Q,
x([0,1[)\subset\Omega_E, x(1)\in\partial\Omega_E\}.
\end{equation}
Then, by Lemma \ref{thm:lem8.2}, $f_{0,1}$ has at least one
minimizer on $X_Q$.
\end{rem}

Using a simple argument, we also have:
\begin{lem}\label{thm:lem84}
The map $d_E:\Omega_E\to\left[0,+\infty\right[$ defined in the
statement of Lemma~\ref{thm:lem8.2} is continuous, and it admits a
continuous extension to $\overline{\Omega}_E$ by setting $d_E=0$
on $\partial\Omega_E$.\qed
\end{lem}

Now we shall study the map
\begin{equation}\label{eq:psi}
\psi(y)=\frac12 d_E^2(y),
\end{equation}
proving that it is $C^2$ and satisfies a convex condition when $y$
is nearby $\partial\Omega_E$.

\begin{prop}\label{prop:uniq}
If $Q$ is sufficiently close to $\partial\Omega_E$ then the
minimizer of the functional \eqref{eq:funct01} in the space $X_Q$
is unique.
\end{prop}

\begin{proof}
Let $z=z(t,0,Q)$ the solution of the Cauchy problem
\begin{equation}\label{eq:ode}
\left\{
\begin{aligned}
&\dot z(t)=J\cdot\,D_z H(z(t))\\
&z(0)=(0,Q),\quad Q\in\partial\Omega_E,
\end{aligned}
\right.
\end{equation}
where $H$ is the Hamiltonian function \eqref{eq:hamfun}, and $J$
is the matrix
\[
J=
\begin{pmatrix}
0 & -I_m \\ I_m & 0
\end{pmatrix}
\]
and $I_m$ is the $m\times m$ identity matrix. Since $V$ and
$a_{ij}$ are $C^2$, $z=(p,q)$ is of class $C^1$ with respect to
$(t,Q)$, therefore $\dot z=\dot z(t,Q)$ is of class $C^1$ with
respect to $(t,Q)$ so $\dot q=\dot q(t,Q)$ is $C^1$. Since $\dot
q=\dot q(0,Q)=0$, in a neighborhood of a fixed point $Q_0\in\partial\Omega_E$ it
is
\begin{equation}\label{eq:eq3.7}
\dot q(t,Q)=t\ddot q(0,Q_0)+\varphi(t,Q)=-t\nabla
V(Q_0)+\varphi(t,Q)
\end{equation}
where $\varphi$ is of class $C^1$ and $\text d\varphi(0,Q_0)=0$.
Moreover
\begin{equation}\label{eq:eq3.8}
q(t,Q)=Q-\frac{t^2}{2}\nabla V(Q_0)+\varphi_0(t,Q)
\end{equation}
where $\varphi_0(t,Q)=\int_0^t\varphi(s,Q)\,\text ds$. Then, if
$\{y_1,\ldots,y_{m-1}\}$ is a coordinate system of $V^{-1}(E)$ in
a neighborhood of $Q_0$, by \eqref{eq:eq3.8} we deduce that,
setting $\tau=t^2$, the set $\{y_1,\ldots,y_{m-1},\tau\}$ is a
local coordinate system on the manifold with boundary
$\partial\Omega_E$ and $(\tau,Q)\mapsto q(\tau,Q)$ defines a local
chart.

Then, due to the compactness of $\partial\Omega_E$, and denoted by
$\dist(\cdot,\cdot)$ the distance induced by $g$, there exists
$\bar\rho>0$ having the following property:
\begin{equation}\label{eq:property}
\begin{matrix}
\text{\textsl{$\forall y\in\Omega_E$ with
$\dist(y,\partial\Omega_E)\le\bar\rho$ there exists a unique
solution $(p_y,q_y)$ of
\eqref{eq:HS}}} \\
\text{\textsl{ with energy $E$, and a unique $t_y>0$ such that
$q_y(0)\in\partial\Omega_E$, $q_y(t_y)=y$.}}
\end{matrix}
\end{equation}
Then, by Proposition \ref{thm:MauJac}, $\forall y\in\Omega_E$ with
$\dist(y,\partial\Omega_E)\le\bar\rho$ there exists a unique
minimizer $\gamma_y$ for $f_{0,1}$ on $X_y$.
\end{proof}

\begin{rem}\label{rem:rem3.5}
Note that $q_y(t)=q(t,Q_y)$ where $Q_y$ is implicitly defined by
$q(t_y,Q_y)=y$. By the variable change used in Proposition
\ref{thm:MauJac}, it turns out that
\begin{equation}\label{eq:eq3.10}
q(t,Q_y)=\gamma_y(1-\sigma),\qquad\text{where\
}t(\sigma)=\psi(y)\,\int_0^\sigma\frac{1}{E-V(\gamma_y(\tau))}\,\text
d\tau.
\end{equation}
In particular, since $\sigma=\sigma(t)$ is the inverse of
$t(\sigma)$ we have
\begin{equation}\label{eq:eq3.11}
\psi(y)\dot q(t_y,Q_y)=-(E-V(y))\dot\gamma_y(0).
\end{equation}
Note also that $t_y=\sqrt{\tau_y}$ is of class $C^1$ when
$\tau_y>0$ since $(\tau,Q)$ is a local coordinate system.
\end{rem}

In the following result we are assuming
$\overline{\Omega}_E\subset\R^m$.
\begin{prop}\label{prop:prop3.6}
Let $\bar\rho$ satisfy property \eqref{eq:property}. Whenever
$0<\dist(y,\partial\Omega_E)\le\bar\rho$, $\psi$ is differentiable
at $y$ and
\begin{equation}\label{eq:dpsi}
\text
d\psi(y)[\xi]=-\big(E-V(y)\big)g\big(\dot\gamma_y(0),\xi\big)\qquad\forall\xi\in\R^m.
\end{equation}
\end{prop}

\begin{proof}
Given the local nature of the result, it will not be restrictive
to assume that $M$ is topologically embedded as an open subset of $\R^m$.
Consider
\[
v_\xi(s)=(1-2s)^+\xi,
\]
where $(\cdot)^+$ denotes the positive part. For $\varepsilon$
sufficiently small (with respect to $\xi$) the curve
$\gamma_y(s)+\varepsilon v_\xi(s)$ belongs to
$X_{y+\varepsilon\xi}$ (see \eqref{eq:Xspace}). Then, by the
definition of $\psi$ as minimum value,
\[
\psi(y+\varepsilon\xi)\le f_{0,1} (\gamma_y+\varepsilon v_\xi)
\]
and therefore
\[
\psi(y+\varepsilon\xi)-\psi(y)\le f_{0,1}(\gamma_y+\varepsilon
v_\xi)-f_{0,1}(\gamma_y).
\]
Now
\begin{multline*}
\lim\limits_{\varepsilon\to 0} \frac{1}{\varepsilon}\left(
f_{0,1}(\gamma_y+\varepsilon v_\xi)-f_{0,1}(\gamma_y)\right)=\\
\int_0^1\big(E-V(\gamma_y)\big)g\big(\dot\gamma_y,\Ddt
v_\xi\big)-\frac 12 g\big(\nabla
V(\gamma_y),v_\xi\big)g\big(\dot\gamma_y,\dot\gamma_y\big)\,\text
ds
\end{multline*}
uniformly as $|\xi|\le 1$. Moreover, since $v_\xi=0$ in the
interval $[\frac 12,1]$, using the differential equation satisfied
by $\gamma_y$ and integrating by parts gives
\begin{multline*}
\int_0^1\big(E-V(\gamma_y)\big)g\big(\dot\gamma_y,\Ddt
v_\xi\big)-\frac 12 g\big(\nabla
V(\gamma_y),v_\xi\big)g\big(\dot\gamma_y,\dot\gamma_y\big)\,\text
ds=\\
-\big(E-V(\gamma_y(0))\big)g\big(\dot\gamma_y(0),v_\xi(0)\big)=-\big(E-V(y)\big)
g\big(\dot\gamma_y(0),\xi\big).
\end{multline*}
Therefore, uniformly as $|\xi|\le 1$,
\begin{equation}\label{eq:eq3.13}
\lim\sup\limits_{\varepsilon\to 0^+}\frac 1\varepsilon
\left(\psi(y+\varepsilon v_\xi)-\psi(y)\right)+\big(E-V(y)\big)
g(\dot\gamma_y(0),\xi\big)\le 0.
\end{equation}
Moreover, since
$\psi(y+\varepsilon\xi)=f_{0,1}(\gamma_{y+\varepsilon\xi})$ and
$\psi(y)\le f_{0,1}(\gamma_{y+\varepsilon\xi}-\varepsilon v_\xi)$
one has
\begin{multline}\label{eq:eq3.14}
\psi(y+\varepsilon\xi)-\psi(y)\ge
f_{0,1}(\gamma_{y+\varepsilon\xi})-f_{0,1}(\gamma_{y+\varepsilon\xi}-\varepsilon
v_\xi)=\\
\varepsilon \langle
f'_{0,1}(\gamma_{y+\epsilon\xi}),v_\xi\rangle_1-\frac{\varepsilon^2}{2}\langle
f''_{0,1}(\gamma_{y+\varepsilon\xi}-\vartheta_\varepsilon
\varepsilon v_\xi)[v_\xi],v_\xi\rangle_1,
\end{multline}
for some $\vartheta_\varepsilon\in ]0,1[$. Here
$\langle\cdot,\cdot\rangle_1$ denotes the standard scalar product
in $H^1$ and $f'$, $f''$ are respectively gradient and Hessian
with respect to $\langle\cdot,\cdot\rangle_1$.

Now, it is  $\gamma_{y+\varepsilon\xi}(0)=y+\varepsilon_\xi$ and
$y\not\in V^{-1}(E)$. Moreover, by the uniqueness of the minimizer
it is not difficult to prove that, $\forall\delta>0$
$\exists\varepsilon(\delta)>0$ such that
\[
\dist(\gamma_{y+\varepsilon\xi}(s),\gamma_y(s))\le\delta\qquad\text{for
any\ } \varepsilon\in ]0,\varepsilon(\delta)],\,|\xi|\le
1,\,s\in[0,1].
\]
Then, since $\gamma_y$ is uniformly far from $V^{-1}(E)$ on the
interval $[0,\frac 12]$, the same holds for
$\gamma_{y+\varepsilon\xi}$ whenever $\varepsilon$ is small and
$|\xi|\le 1$. Thus, recalling the definition of $d_E$ in Lemma
\ref{thm:lem8.2}, the conservation law satisfied by the minimizer
$\gamma_{y+\varepsilon\xi}$ is
\[
\big(E-V(\gamma_{y+\varepsilon\xi})\big)
g\big(\dot\gamma_{y+\varepsilon\xi},\dot\gamma_{y+\varepsilon\xi}\big)=
d_E^2(y+\varepsilon\xi).
\]
This implies the existence of a constant $C>0$ such that
\[
\int_0^{1/2}g\big(\dot\gamma_{y+\varepsilon\xi},\dot\gamma_{y+\varepsilon\xi}\big)
\,\text ds\le C
\]
for any $\varepsilon$ small and $|\xi|\le 1$.

Therefore $\langle
f''_{0,1}(\gamma_{y+\varepsilon\xi}-\vartheta_\varepsilon\varepsilon
v_\xi)[v_\xi],v_\xi\rangle_1$ is uniformly bounded with respect to
$\varepsilon$ small and $|\xi|\le 1$, due to $v_\xi=0$ on $[\frac
12,1]$, and by \eqref{eq:eq3.14} we get
\begin{equation}\label{eq:eq3.15}
\lim\limits_{\varepsilon\to 0}\frac
1\varepsilon\left(f_{0,1}(\gamma_{y+\varepsilon\xi})-f_{0,1}(\gamma_{y+\varepsilon\xi-\varepsilon
v_\xi})\right)=\lim\limits_{\varepsilon\to 0}\langle
f'_{0,1}(\gamma_{y+\varepsilon\xi}),v_\xi\rangle_1
\end{equation}
uniformly as $|\xi|\le 1$.

Now, using the differential equation \eqref{eq:pr1} satisfied by
$\gamma_{y+\varepsilon\xi}$ and integrating by parts one obtains
\[
\langle f'_{0,1}(\gamma_{y+\varepsilon\xi}),v_\xi\rangle_1=
-\big(E-V(y+\varepsilon\xi)\big)g\big(\dot\gamma_{y+\varepsilon\xi}(0),\xi\big),
\]
while by \eqref{eq:eq3.11} and the continuity of $\dot q(t_y,Q_y)$
and $\psi(y)$ we have
\begin{equation}\label{eq:eq3.16}
\lim\limits_{\varepsilon\to
0}\big(E-V(y+\varepsilon\xi)\big)\dot\gamma_{y+\varepsilon\xi}(0)=
\big(E-V(y)\big)\dot\gamma_y(0)
\end{equation}
uniformly as $|\xi|\le 1$. Therefore, by
\eqref{eq:eq3.14}--\eqref{eq:eq3.16} it is
\begin{equation}\label{eq:eq3.17}
\lim\inf\limits_{\varepsilon\to 0}\frac
1\varepsilon\left(\psi(y+\varepsilon\xi)-\psi(y)\right)+
\big(E-V(y)\big)g\big(\dot\gamma_y(0),\xi\big)\ge 0
\end{equation}
uniformly as $|\xi|\le 1$. Finally, combining \eqref{eq:eq3.13}
and \eqref{eq:eq3.17} one has \eqref{eq:dpsi}.
\end{proof}

\begin{rem}\label{rem:rem3.7}
By \eqref{eq:eq3.11} we deduce that $(E-V(y))\dot\gamma_y(0)$ is
continuous, therefore by \eqref{eq:dpsi}, $\psi$ is of class
$C^1$. Again by \eqref{eq:eq3.11} and the $C^1$--regularity of
$\dot q_y(t_y,Q_y)$ we deduce that $(E-V(y))\dot\gamma_y(0)$ is of
class $C^1$ whenever $y\not\in V^{-1}(E)$, and by \eqref{eq:dpsi}
it turns out that $\psi$ is of class $C^2$.
\end{rem}

In the following proposition we will show that $\psi$ satisfies a
strongly convex assumption nearby $V^{-1}(E)$.
\begin{prop}\label{prop:prop3.8}
There exists $\widehat\rho\le\bar\rho$ with the property that, for
any $y\in\Omega_E$ such that $0<\dist(y,V^{-1}(E))\le\widehat\rho$
the Hessian (with respect to the Jacobi metric $g_E$) of $\Psi$ at $y$ satisfies
\begin{equation}\label{eq:psiconvex}
\mathrm H^\psi(y)[v,v]>0\qquad\forall v\,:\,\text
d\psi(y)[v]=0,\quad v\not=0.
\end{equation}
\end{prop}

\begin{proof}
Recall that
\[
\mathrm H^\psi(y)[v,v]=\frac{\partial^2}{\partial
s^2}\left(\psi(\eta(s))\right)_{\vert s=0},
\]
where $\eta(s)$ is a geodesic with respect to the Jacobi metric
$g_E$, namely a solution of the differential equation
\eqref{eq:pr1} satisfying the initial data conditions
\[
\eta(0)=y,\qquad\dot\eta(0)=\xi.
\]
Now, by \eqref{eq:eq3.11} and \eqref{eq:dpsi}
\[
\text d\psi(\eta(s))[\dot\eta(s)]=-\big(E-V(\eta(s))\big)
g\big(\dot\gamma_{\eta(s)}(0),\dot\eta(s)\big)=
\psi(\eta(s))g\big(\dot
q(t_{\eta(s)},Q_{\eta(s)}),\dot\eta(s)\big).
\]
Since $\lim\limits_{s\to 0}Q_{\eta(s)}=Q_y$, using
\eqref{eq:eq3.7} we can write
\[
\dot q(t,Q_{\eta(s)})=-t\nabla V(y)+\varphi(t,Q_{\eta(s)})
\]
as $\text d\varphi(0,Q_y)=0$, and
\begin{multline*}
\frac{\partial^2}{\partial s^2}\left(\psi(\eta(s))\right)=\\
\psi(\eta(s))\left(g\big(\dot
q(t_{\eta(s)},Q_{\eta(s)}),\dot\eta(s)\big)\right)^2+
\psi(\eta(s))g\big(\dot
q(t_{\eta(s)},Q_{\eta(s)}),\Dds\dot\eta(s)\big)+
\\
\psi(\eta(s)) g\big(-\text dt_{\eta(s)}[\dot\eta(s)]\,\nabla V(y)+
\frac{\partial\varphi}{\partial t}(t_y,Q_{\eta(s)}) \text
dt_{\eta(s)}[\dot\eta(s)]+\frac{\partial\varphi}{\partial
Q}\,\frac{\partial Q}{\partial\eta}[\dot\eta(s)],\dot\eta(s)\big).
\end{multline*}
Since $\eta(s)$ satisfies \eqref{eq:pr1} and $\text
d\varphi(0,Q_y)=0$, it suffices to show that for any $y$
sufficiently close to $\partial\Omega$,
\begin{multline*}
\psi(\eta(s))\left(g\big(\dot q(t_y,Q_y),v\big)\right)^2+\psi(y)\text dt_y[v]
g\big(-\nabla V(y),v\big)+\\
\frac{\psi(y)}{E-V(y)}\left( g\big(\nabla V(y),v\big) g\big(\dot
q(t_y,Q_y),v\big) -\frac 12 g\big(\dot q(t_y,Q_y),\nabla V(y)\big)
g(v,v)\right)>0
\end{multline*}
for any $v$ such that $\text d\psi(y)[v]=0$. This means that
$g\big(\dot q(t_y,Q_y),v\big)=0$ so it will suffice to show
\begin{equation}\label{eq:eq3.19}
\sup\limits_{|v|=1}|\text dt_y[v]| g\big(\nabla V(y),\nabla
V(y)\big)^{1/2} -\frac 1{2(E-V(y))}g\big(\dot q(t_y,Q_y),\nabla
V(y)\big)>0
\end{equation}
for any $y$ close to $V^{-1}(E)$.

Since $q(t_y,Q_y)=y$ we get
\[
\text dt_y[v]\dot q(t_y,Q_y)+\frac{\partial q}{\partial
Q}\,\frac{\partial Q_y}{\partial y}[v]=v.
\]
Moreover, $\frac{\partial q}{\partial Q}(t_y,Q_y)$ goes to the
identity map as $y$ tends to $\partial\Omega$, while
$\frac{\partial Q_y}{\partial y}[v]$ tends to $v$ uniformly as
$|v|\le 1$, since $(0,Q)$ is a coordinate system for $V^{-1}(E)$.
Then, as $y\to V^{-1}(E)$, $\text dt_y[v]\dot q(t_y,Q_y)\to 0$
uniformly in $v$ .

Note that $\frac 12 g(\dot q,\dot q)=E-V(q)$, therefore
\begin{equation}\label{eq:eq3.20}
g\big(\dot q(t_y,Q_y),\dot q(t_y,Q_y)\big)=2\big(E-V(y)\big)
\end{equation}
so
\begin{equation}\label{eq:eq3.21}
\lim\limits_{y\to\partial\Omega}\sqrt{E-V(y)}\left|\text
dt_y[v]\right|=0
\end{equation}
uniformly in $|v|\le 1$.

Finally, by \eqref{eq:eq3.7} we have
\[
\lim\limits_{y\to V^{-1}(E)}g\left(\frac{\dot q(t_y,Q_y)}
{\sqrt{g\big(\dot q(t_y,Q_y),\dot q(t_y,Q_y)\big)}}, \frac{\nabla
V(y)}{\sqrt{g\big(\nabla V(y),\nabla V(y)\big)}}\right)=-1
\]
therefore by \eqref{eq:eq3.20}
\begin{equation}\label{eq:eq3.22}
\lim\inf\limits_{y\to V^{-1}(E)} \frac{-g\big(\dot
q(t_y,Q_y),\nabla V(y)\big)}{\sqrt{E-V(y)}}>0
\end{equation}
and combining \eqref{eq:eq3.21} with \eqref{eq:eq3.22} one obtains
\eqref{eq:eq3.19} and the proof is complete.
\end{proof}

By Proposition \ref{prop:prop3.6}, Remark \ref{rem:rem3.7} and
Proposition \ref{prop:prop3.8} one immediately obtains the
following proposition, which is the main result of the section:
\begin{teo}\label{thm:lem8.5}
Let $E$ be a regular value for $V(x)$, and
let $d_E:\Omega\to\left[0,+\infty\right[$ be the map defined in
the statement of Lemma~\ref{thm:lem8.2}, and assume that $\overline \Omega_E$ is compact.
There exists a positive
number $\delta_*$ such that, setting:
\[\Omega_*=\big\{x\in\Omega_E\,:\,d_E(x)>\delta_*\big\},\]
the following statements hold:
\begin{enumerate}
\item $\partial\Omega_*$ is of class $C^2$; \item
$\overline\Omega_*$ is omeomorphic to $\overline\Omega_E$; \item
$\overline\Omega_*$ is {\em strongly\/} concave relatively to the
Jacobi metric $g_E$; \item if $x:[0,1]\to\overline\Omega_*$ is an
orthogonal geodesic chord in $\overline\Omega_*$ relatively to the
Jacobi metric $g_E$, then there exists
$[\alpha,\beta]\supset[0,1]$ and a unique extension $\widehat
x:[\alpha,\beta]\to\overline\Omega$ of $x$ with $\widehat x\in
H^1\big([\alpha,\beta],\overline\Omega\big)$ satisfying:
\begin{itemize}
\item assumption \eqref{eq:i} of Proposition~\ref{thm:MauJac} on
the interval $[\alpha,\beta]$; \item $\widehat x(s)\in
d_E^{-1}\big(\left]-\delta_*,0\right[\big)$ for all
$s\in\left]\alpha,0\right[\bigcup\left]1,\beta\right[$; \item
$V\big(\widehat x(\alpha)\big)= V\big(\widehat x(\beta)\big)=E$.
\end{itemize}
\end{enumerate}
\end{teo}

\begin{rem}\label{rem:rem3.10}
Theorem~\ref{thm:lem8.5} tells us that the study of multiple
brake orbits can be reduced to the study of multiple orthogonal geodesic
chords in a Riemannian manifold with regular and strongly concave boundary.
\end{rem}

\subsection{The Jacobi distance near a nondegenerate maximum point of the potential.}

Let us now assume that $x_0\in M$ is a nondegenerate maximum point
of $V$, with $V(x_0)=E$, and let us make the following assumptions:
\begin{itemize}
 \item $V^{-1}\big(\left]-\infty,E\right]\big)$ is compact;
 \item $V^{-1}(E)\setminus\{x_0\}$ is a regular embedded hypersurface
 of $M$.
 \end{itemize}
 We will show how to get rid of the singularity of the Jacobi metric
 at $x_0$, while the singularity on $V^{-1}(E)\setminus\{x_0\}$
can be removed as in the case of brake orbits, 
using Theorem~\ref{thm:lem8.5}.

First, we need a preparatory result. Let $\delta>0$ be fixed in
such a way that the set:
\[\Big\{p\in M:V(p)>E-\delta\Big\}\]
has precisely two connected components; let $\Omega_\delta$ denote
the connected component of the point $x_0$.
\begin{lem}\label{thm:lem9.2}
Let $Q\in\Omega_\delta\setminus\{x_0\}$ be fixed; then, the
infimum:
\begin{multline}\label{eq:infimum}
d_E(Q):=\inf\Big\{\left[\int_0^1(E-V(x))g(\dot x ,\dot x)\,\text dt\right]^{1/2}:\\
x\in C^0\big([0,1],\overline{\Omega_\delta}\big)\cap
H^1_{\textrm{loc}}\big(\left[0,1\right[,\overline{\Omega_\delta}\big),
\ x(0)=Q,\ x(1)=x_0\Big\}
\end{multline}
is attained on some curve $\gamma_Q$ with the property
$(E-V(\gamma_Q)) g(\dot\gamma_Q,\dot\gamma_Q)$ constant and
$\gamma_Q([0,1[)\subset\overline{\Omega_\delta}\setminus\{x_0\}$.
Moreover
\begin{eqnarray}\label{eq:eq9.3}
&&\lim_{Q\to x_0}d_E(Q)=0,\\
\label{eq:eq9.4}
&&\lim_{Q\to x_0}\left[\sup_{s\in[0,1]}\dist\big(\gamma_Q(s),x_0\big)\right]=0,\label{eq:eq9.5} \\
\end{eqnarray}
In particular, for $Q$ sufficiently close to  $x_0$,
\begin{equation}\label{eq:eq9.7}
\gamma_Q\big(\left[0,1]\right)\subset\Omega_\delta,
\end{equation}
so it is of class $C^2$ and satisfies assumption \eqref{eq:2.4bis}
of Proposition \ref{thm:MauJacnearmax} on the interval
$[a,b]=[0,1]$.
\end{lem}
\begin{proof}
Let $x_n\in C^0\big([0,1],\overline{\Omega_\delta}\big)\cap
H^1\big(\left[0,1\right[,\overline{\Omega_\delta}\big)$ be a
minimizing sequence for the length functional
$\int_0^1\left[(E-V(x))g(\dot x ,\dot x)\right]^{1/2}\,\text dt$,
leaving $(E-V(x))g(\dot x ,\dot x)$ constant. Choose $\rho>0$ such
that $\dist(Q,x_0)>\rho$ and, for all $n\in\N$, define
$\alpha_1^n\in\left]0,1\right[$ to be the first instant $s$ such
that $\dist\big(x_n(s),x_0\big)=\rho$.

The sequence $\alpha_1^n$ stays away from $0$ and $1$, because for
all interval $I\subset x_n^{-1}\big([\frac\rho2,\rho]\big)$ the
integral $\int_Ig(\dot x_n,\dot x_n)\,\mathrm ds$ is bounded. We
can therefore find a subsequence $\alpha_1^{n_k}$ converging to
$\alpha_1\in\left]0,1\right[$. Furthermore, since
$\int_0^{\alpha_1}g(\dot x_n,\dot x_n)\,\mathrm ds$ is bounded,
taking a subsequences $x_n^1$ we can assume that $x_{n}^1$ is
$H^1$-weakly and uniformly convergent to some $x_1\in
H^1\big([0,\alpha_1],\overline{\Omega_\delta})$; then,
$\dist\big(x(\alpha_1),x_0\big)=\rho$. Repeating the construction,
we can find $\alpha_2\in\left]\alpha_1,1\right[$ and a subsequence
$x_n^2$ of $x_n^1$ which is $H^1$-weakly and uniformly convergent
to a curve $x_2\in H^1
\big([0,\alpha_2],\overline{\Omega_\delta}\big)$ with
$\dist\big(x(\alpha_2),x_0\big)=\frac\rho2$ and $x_{2\vert
[0,\alpha_1]}=x_1$. Iteration of this construction yields a
weak-$H^1$ limit of $x_n^n$, which is a curve $x\in
H^1_{\textrm{loc}}\big(\left[0,\bar\alpha\right[,\overline{\Omega_\delta}\big)$,
where $\bar\alpha=\lim\limits_k\alpha_k$, and
$\dist\big(x(\alpha_k),x_0\big)=\frac\rho{2^k}$.

Now, for all $k\ge1$:
\begin{multline*}
\int_0^{\alpha_k}\left(\big(E-V(x)\big)g(\dot x,\dot
x)\right)^{1/2}\,\mathrm ds
\le\liminf_{n\to\infty}\int_0^{\alpha_k}\left(\big(E-V(x_n)\big)g(\dot
x_n,\dot x_n)\right)^{1/2}\, \mathrm ds\\ \le
\liminf_{n\to\infty}\int_0^{1}\left(\big(E-V(x_n)\big)g(\dot
x_n,\dot x_n)\right)^{1/2}\, \mathrm ds=d_E(Q),
\end{multline*}
hence:
\[
\int_0^{\bar\alpha}\left(\big(E-V(x)\big)g(\dot x,\dot
x)\right)^{1/2}\,\mathrm ds=
\lim_{k\to\infty}\int_0^{\alpha_k}\left(\big(E-V(x)\big)g(\dot
x,\dot x)\right)^{1/2}\,\mathrm ds \le d_E(Q).
\]
and we can assume, as usual, $(E-V(x))g(\dot x ,\dot x)$ constant
(and positive since $Q\not= x_0$). The curve $x$ can be extended
continuously to $\overline\alpha$ by setting
$x(\overline\alpha)=x_0$. Indeed, if by contradiction there exists
a sequence $\beta_n<\alpha_n<\overline\alpha$ such that
$\lim_k\beta_k=\overline\alpha$ and a positive number
$\overline\nu$ such that $\dist(x(\beta_k),x_0)\ge\overline\nu$,
there exist $\beta_k^1\in ]\beta_k,\alpha_k[$ such that
$\dist(x(\beta_k^1),x_0)=\frac{\overline\nu}2$ and
$\dist(x(s),x_0)\ge\frac{\overline\nu}2$, $\forall
s\in[\beta_k^1,\beta_k]$. But $E-V(x(s))$ is far from zero in
$[\beta_k^1,\beta_k]$ therefore $g(\dot x,\dot x)\le K\in\R^+$ on
$[\beta_k^1,\beta_k]$ for some $K$, and then
\[
\frac{\overline\nu}2\le\dist(x(\beta_k^1),x(\beta_k))\le\int_{\beta_k^1}^{\beta_k}
g(\dot x,\dot x)\,\text dt\le K(\beta_k-\beta_k^1)\longrightarrow
0
\]
which is a contradiction.

Clearly, up to reparameterizations on $x$ we can assume
$\overline\alpha=1$ and
$x([0,1[)\subset\overline{\Omega_\delta}\setminus\{x_0\}$. Taking
$\gamma_Q=x$ we have the existence of a minimizer satisfying the
conservation law $(E-V(\gamma_Q)) g(\dot\gamma_Q,\dot\gamma_Q)$
constant.

Now, taking a chord $C_Q$ joining $Q$ and $x_0$ we have that
$l(C_Q)\to 0$ as $Q\to x_0$, and since $d_E(Q)\le l(C_Q)$ we
obtain \eqref{eq:eq9.3}.

Moreover, if by contradiction \eqref{eq:eq9.5} does not hold for
any $Q$ sufficiently close to $x_0$, there exists $s_Q$ such that
\[
\dist(\gamma_Q(s_Q),x_0)\ge\overline\nu>0.
\]
Let $t_Q>s_Q$ such that
$\dist(\gamma_Q(t_Q),x_0)=\frac{\overline\nu}2$ and
$\dist(\gamma_Q(s),x_0)\ge\frac{\overline\nu}2\,\forall s\in
[s_Q,t_Q]$. Since $g(\dot\gamma_Q,\dot\gamma_Q)$ is bounded in
$[s_Q,t_Q]$ it must be $t_Q-s_Q$ far from zero as $Q\to x_0$. But
also $E-V(\gamma_Q)$ and $g(\dot\gamma_Q,\dot\gamma_Q)$ are far
from zero in $[s_Q,t_Q]$ so we deduce that
\[
\int_{s_Q}^{t_Q} \left(\int_0^1(E-V(x))g(\dot\gamma_Q
,\dot\gamma_Q)\,\text dt\right)^{1/2} \text{\ far from zero}
\]
which is in contradiction with \eqref{eq:eq9.3}.

Note that \eqref{eq:eq9.5} immediately implies \eqref{eq:eq9.7}
and since $\gamma_Q$ is a minimizer satisfying $(E-V(\gamma_Q))
g(\dot\gamma_Q,\dot\gamma_Q)$ constant, we immediately see that
\eqref{eq:2.4bis} is satisfied in the interval $[0,1]$.
\end{proof}

As for Lemma \ref{thm:lem84} a simple argument shows

\begin{lem}\label{lem:lem3.11}
The map $d_E:\Omega_\delta\to [0,+\infty[$ defined in the
statement of Lemma \ref{thm:lem9.2} is continuous.
\end{lem}

For any $y$ sufficiently close to $x_0$, let $q_y$ be the
reparameterization of $\gamma_y$ given by Proposition
\ref{thm:MauJacnearmax}. We have
\begin{equation}\label{eq:eq3.26}
\left\{
\begin{aligned}
&\Dds\dot q_y+\nabla Y(q_y)=0\\
&q_y(0)=y\\
&\lim_{t\to+\infty}q_y(t)=x_0\\
&\lim_{t\to+\infty}\dot q_y(t)=0.
\end{aligned}
\right.
\end{equation}
The following estimate holds

\begin{prop}\label{prop:est1}
Let $q_y$ be as above. Then there exists $\bar\rho$ and a constant
$\alpha>0$ such that
\begin{equation}\label{eq:eq3.27}
\dist(q_y(t),x_0)\le\dist(y,x_0)e^{-\alpha t}
\end{equation}
for any $y$ such that $\dist(y,x_0)\le\bar\rho$.
\end{prop}

To obtain the above result we need the following maximum principle
in $\R$.

\begin{lem}\label{lem:max}
Let $\varphi:[0,+\infty[\to\R$ be a $C^2$ map with
$\lim_{t\to+\infty}\varphi(t)=0$. Let $\nu>0$ such that
$\varphi''(t)\ge\nu\varphi(t),\,\forall t\ge 0$. Then
$\varphi\le\varphi(0) e^{-\sqrt\nu t}$.
\end{lem}

\begin{proof}
Consider the map $\psi=\varphi-\varphi_0$ where
$\varphi_0(t)=\varphi(0)e^{-\sqrt\nu t}$. Clearly
$\psi(0)=\lim_{t\to+\infty}\psi(t)=0$ and so $\psi$ has a global
maximum at some $\bar t\in [0,+\infty[$. If $\bar t>0$ then
$\psi(\bar t)\le\frac 1\nu\psi''(\bar t)\le 0$.
\end{proof}

\begin{rem}\label{rem:rev}
Clearly, an analogous result as in the above Lemma \ref{lem:max}
holds, reversing all inequalities.
\end{rem}

\begin{proof}[Proof of Proposition \ref{prop:est1}]
Let $q$ be a solution of \eqref{eq:eq3.26} (with $q(0)=y$), and
let $\varphi(t)=\frac 12\dist(q(t),x_0)^2$. By \eqref{eq:eq9.5} we
can choose $\bar\rho$ sufficiently small so that
\[
\dist(q(t),x_0)<\rho_0,\quad\text{\ for any } t\ge 0,
\]
where $\rho_0$ is chosen so that the function $\mathfrak
d(z)=\frac12\dist(z,x_0)^2$, in the open ball $B(x_0,\rho_0)$ of
center $x_0$ and radius $\rho_0$, is of class $C^2$, strictly
convex and, called $x_z$ the unique minimal geodesic with respect
to $g$ such that $x_z(0)=x_0$, $x_z(1)=z$ (see \cite{docarmo}),
one has
\[
\nabla\mathfrak d(z)=\dot x_z(1).
\]
Now $\varphi'(t)=g\big(\nabla\mathfrak d(q(t)),\dot q(t)\big)$ and
\[
\varphi''(t)=\mathrm H^{\mathfrak d}(q(t))[\dot q(t),\dot
q(t)]+g\big(\nabla\mathfrak d(q(t)),\Ddt\dot q(t)\big)\ge
g\big(\nabla\mathfrak d(q(t)),\nabla V(q(t))\big).
\]
Now, take $z$ in $B(x_0,\rho_0)$, consider the minimal geodesic
$x_z$ as above, and define the map
\[
\rho(s):=g\big(\nabla\mathfrak d(x_z(s)),-\nabla V(x_z(s))\big).
\]
By the choice of $x_z$ it is $\nabla\mathfrak d(x_z(s))=s\,\dot
x_z(s)$, so
\begin{multline*}
\dot\rho(s)=g\big(\dot x_z(s),-\nabla V(x_z(s))\big)
-s\,\mathrm H^V(x_z(s))[\dot x_z(s),\dot x_z(s)]\ge \\
g\big(\dot x_z(s),-\nabla V(x_z(s))\big)+s\nu g(\dot x_z(s),\dot
x_z(s)\big)
\end{multline*}
for a suitable choice of $\nu$ ($x_0$ is a nondegenerate maximum
point). Since $\rho(0)=0$ then
\begin{multline*}
g\big(\nabla\mathfrak d(z),-\nabla
V(z)\big)=\varphi(1)=\int_0^1\dot\rho(s)\,\text ds\ge \\
\int_0^1 g\big(\dot x_z(s),-\nabla V(x_z(s))\big)+s\nu g\big(\dot
x_z(s),\dot x_z(s)\big)\,\text ds=\\
-V(x_z(s)\big\vert_{s=0}^{s=1}+\nu\dist(z,x_0)^2\int_0^1 s\,\text
ds=\\
(E-V(z))+\frac\nu 2\dist(z,x_0)^2\ge\frac\nu 2\dist(z,x_0)^2,
\end{multline*}
where $V(x_0)=E$ has also been used. Therefore
$\varphi''(t)\ge\frac\nu 2\dist(q(t),x_0)^2=\nu\,q(t)$, and by
Lemma \ref{lem:max}
\[
\dist(q(t),x_0)^2\le \dist(q(0),x_0)e^{-\sqrt\nu t},
\]
and \eqref{eq:eq3.27} follows taking the square root of both
members above.
\end{proof}

The regularity of the distance function from $x_0$ with respect to
the Jacobi metric is based on the following proposition.

\begin{prop}\label{prop:prop3.14}
For any $y$ close to $x_0$ there exists a unique $q_y$ satisfying
\eqref{eq:eq3.26}. Moreover, the map
\begin{equation}\label{eq:eq3.28}
q\longmapsto\dot q_y(0)
\end{equation}
is of class $C^1$ and its differential satisfies $\text d\dot
q_y(0)[v]=\dot\xi(0)$, where $\xi(t)$ is the unique solution of
\begin{equation}\label{eq:eq3.29}
\left\{
\begin{aligned}
&\Ddtt\xi(t)+R(\dot q_y,\xi(t))\dot q_y+\mathrm L^V(q_y)\xi(t)=0\\
&\xi(0)=0\\
&\lim_{t\to+\infty}\xi(t)=\lim_{t\to+\infty}\dot\xi(t)=0
\end{aligned}
\right.
\end{equation}
where $\Ddtt\xi$ is the second covariant derivative and
$R(\cdot,\cdot)$ the Riemann tensor with respect to $g$, and
$\mathrm L^V(x)[v]\in T_x M$ is the vector defined through
$g\big(\mathrm L^V(x)[v],w)=H^V(x)[v,w]$ for all $w\in T_x M$.
\end{prop}

\begin{proof}
Consider the ball $B(x_0,\rho)$, with $\rho>0$ small, and the
spaces
\begin{equation*}
X_2=\{q\in C^2(\R^+,\overline{B(x_0,\rho)})\,:\,\lim_{t\to+\infty}
q(t)=x_0,\,\lim_{t\to+\infty} \dot q(t)=\lim_{t\to+\infty} \ddot
q(t)=0\}
\end{equation*}
with the norm (we can assume to work in a local chart)
\begin{equation}\label{eq:eq3.29bis}
\Vert
q_2-q_1\Vert:=\sup_{t\in\R^+}|q_2(t)-q_1(t)|+\sup_{t\in\R^+}|\dot
q_2(t)-\dot q_1(t)|+ \sup_{t\in\R^+}|\ddot q_2(t)-\ddot q_1(t)|
\end{equation}
and
\begin{equation*}
X_0=\{q\in C^0(\R^+,\R^m)\,:\,\lim_{t\to+\infty}|q(t)|=0\}
\end{equation*}
with the norm
\[
\Vert q_2-q_1\Vert:=\sup_{t\in\R^+}|q_2(t)-q_1(t)|,
\]
that are clearly Banach spaces. Now, consider the open  set
\[
A_2=\{q\in X_2\,:\,\sup_{t\in\R^+}\dist(q(t),x_0)<\rho\}\subset
X_2
\]
and the map
\[
F:A_2\times B(x_0,\rho)\longrightarrow X_0\times\R^m
\]
given by
\[
F(q,y)=(\Ddt\dot q+\nabla V(q),q(0)-y).
\]
Thanks to the behaviour at infinity, we can use the same standard
arguments exploited in finite intervals to prove that $F$ is
differentiable and (see \cite{docarmo})
\[
\text dF(q,y)[\xi,v]=\big(\Ddtt\xi+R(\dot q,\xi)\dot q+\mathrm
L^V(q)[\xi],\xi(0)-v\big).
\]
Moreover, thank again to the behaviour at infinity, it is a
straight check to verify that $\text dF(q,y)$ is continuous
(recall that $g$ and $V$ are of class $C^2$).

Now consider $\frac{\partial F}{\partial
q}(x_0,0)[\xi]=\big(\ddot\xi+\mathrm L^V(0)[\xi],\xi(0)\big)$
where $x_0$ denotes the constant curve with image $x_0$. We claim
that
\begin{equation}\label{eq:eq3.30}
\frac{\partial F}{\partial q}(x_0,0):X_2\longmapsto X_0\times\R^m
\end{equation}
is an isomorphism.

Recalling the definition of $\mathrm L^V$, and since $\mathrm
H^V(0)$ is symmetric and negative definite, using a base
consisting of eigenvectors for $\mathrm H^V(0)$, it is sufficient
to show that for any function $h\in C^0(\R^+,\R)$ such that
$\lim_{t\to+\infty} h(t)=0$ and for any $\theta\in\R$, the
solution of
\begin{equation}\label{eq:eq3.31}
\left\{
\begin{aligned}
&\ddot x-\alpha^2 x=h\\
&x(0)=\theta\\
&\lim_{t\to+\infty}x(t)=\lim_{t\to+\infty}\dot x(t)=0
\end{aligned}
\right.
\end{equation}
exists and is unique (where $x:\R^+\to\R$).

The general solution of the differential equation above is
\[
x(t)=\left(a+\frac{1}{2\alpha}\int_0^t h(s)e^{-\alpha s}\,\text
ds\right)e^{\alpha t}+\left(b-\frac{1}{2\alpha}\int_0^t h(s)
e^{\alpha s}\,\text ds\right)e^{-\alpha t}.
\]
Since $\lim\limits_{t\to+\infty} h(t)=0$ it is
\[
\lim\limits_{t\to+\infty}\frac 1{2\alpha}\left(\int_0^t h(s)
e^{\alpha s}\,\text ds\right) e^{-\alpha t}=0
\]
then $\lim\limits_{t\to+\infty}x(t)=0$ only if we choose
\[
a=-\frac 1{2\alpha}\int_0^{+\infty}h(s) e^{-\alpha s}\,\text ds.
\]
With such a choice indeed $\lim\limits_{t\to+\infty}
x(t)=\lim\limits_{t\to+\infty}\dot x(t)=0$, while $x(0)=\theta$
for
\[
b=\theta-a=\theta+\frac 1{2\alpha}\int_0^{+\infty}h(s)e^{-\alpha
s}\,\text ds,
\]
proving that the solution of \eqref{eq:eq3.31} exists and is
unique, and 
therefore the map defined in \eqref{eq:eq3.30} is an isomophism.

Then, by the Implicit Function Theorem and Proposition
\ref{prop:est1} we have the uniqueness of $q_y$ for any $y$ close
to $x_0$ and its $C^1$--differentiability in $X_2$. In particular
the map \eqref{eq:eq3.28} is of class $C^2$. Denoting by $\xi$ the
differential $\text dq_y[v]$, and differentiating the expression
$F(q_y,y)\equiv 0$, in particular we obtain that $\xi$ solves
\eqref{eq:eq3.29}. Since, has we have already seen, the solution
exists and is unique for $y=x_0$, Proposition \ref{prop:est1}
ensures that this remains true for $y$ close to $x_0$ also.

Finally, $C^1$--regularity of $q_y$ with respect to the norm
\eqref{eq:eq3.29bis} immediately implies that
\[
\text d\dot q_y(0)[v]=\dot\xi(0),
\]
where $\xi$ is the solution of \eqref{eq:eq3.29}, and then $\text
dq_y[v](t)=\xi(t)$.
\end{proof}

Now set
\begin{equation}\label{eq:eq3.32}
\psi(y)=\frac 12 d_E(y)^2
\end{equation}
where $l$ is the map defined in \eqref{eq:infimum} of Lemma
\ref{thm:lem9.2}. Thanks to the above proposition we can repeat
the proof of Proposition \ref{prop:prop3.6} to get its counterpart in
the case of a nondegenerate maximum point.

\begin{prop}\label{prop:prop3.15}
There exists $\bar\rho>0$ such that for any $y$ with
$\dist(y,x_0)\le\bar\rho$ the map $\psi$ defined in
\eqref{eq:eq3.32} is of class $C^2$ and its differential is given
by
\begin{equation}\label{eq:dpsi3}
\text d\psi(y)[v]
=-(E-V(y))g\big(\dot\gamma_y(0),v\big)=-\psi(y)g\big(\dot
q_y(0),v\big).
\end{equation}
\end{prop}

Note that the variable change used in the proof of Proposition
\ref{thm:MauJacnearmax} yields $q_y(t)=\gamma_y(\sigma)$ where
$t(\sigma)=\psi(y)\int_0^\sigma\frac 1{E-V(\gamma_y(\tau))}\,\text
d\tau$. 

We now are going to show the counterpart of Proposition \ref{prop:prop3.8}. We cannot repeat, of course, the same argument as before: indeed, since $E$ is not a regular value for the potential $V(x)$, the curve $q_y(t)=q(t,Q_y)$ (see Remark \ref{rem:rem3.5}) does not reach the boundary $\partial\Omega$ in a finite amount of time and therefore it cannot be reparameterized in a bounded interval.

\begin{prop}\label{prop:prop3.16}
There exists  $\widehat\rho\le\bar\rho$ such that for any $y$ with
$\dist(y,x_0)\le\widehat\rho$ it is
\[
\mathrm H^\psi(y)[v,v]>0,\qquad\forall v\,:\,\text d\psi(y)[v]=0.
\]
\end{prop}

\begin{proof}
We need to evaluate
\[
\frac{\partial^2}{\partial s^2}\big(\psi(\eta(s))\big)_{\vert
s=0},
\]
where $\eta(s)$ is the geodesic with respect to the Jacobi metric
$g_E$ such that $\eta(0)=y$, $\dot\eta(0)=v$, where $\text
d\psi(y)[v]=0$. We also recall that $\eta(s)$ satisfies equation
\eqref{eq:pr1}. By \eqref{eq:dpsi3}
\begin{multline*}
\frac{\partial^2}{\partial
s^2}\big(\psi(\eta(s))\big)=\frac{\partial}{\partial s}\big(\text
d\psi(\eta(s))[\dot\eta(s)]\big)= \frac\partial{\partial
s}\left(-\psi(\eta(s))g\big(\dot
q_{\eta(s)}(0),\dot\eta(s)\big)\right)=\\
-\text d\psi(\eta(s))[\dot\eta(s)]g\big(\dot
q_{\eta(s)}(0),\dot\eta(s)\big)-\psi(\eta(s))g\big(\Dds(\dot q_{\eta(s)}(0)),\dot\eta(s)\big)-\\
\psi(\eta(s))g\big(\dot q_{\eta(s)}(0),\Dds\dot\eta(s)\big),
\end{multline*}
then, using again \eqref{eq:dpsi3}, and exploiting \eqref{eq:pr1},
one gets
\begin{multline*}
\mathrm H^\psi(y)[v,v]=\psi(y) g\big(\dot
q_y(0),v\big)^2-\psi(y)g\big(\text d\dot q_y(0)[v],v\big)-\\
\frac{\psi(y)}{E-V(y)}\left(-\frac12 g\big(v,v\big)g\big(\dot
q_y(0),\nabla V(y)\big)+g\big(\nabla V(y),v\big) g\big(\dot
q_y(0),v\big)\right).
\end{multline*}
Since $g(\dot q_y(0),v)=\text d\psi(y)[v]=0$, it suffices to show
the existence of $\nu_0>0$ such that
\begin{equation}\label{eq:eq3.33}
\inf\limits_{|v|=1} g\big(\text d\dot
q_y(0)[v],v\big)+\frac{g\big(\dot q_y(0),\nabla
V(y)\big)}{2\big(E-V(y)\big)}\ge\nu_0
\end{equation}
for any $y$ close sufficiently to $x_0$. Let us consider the map
$\mu(t)=g\big(\dot q_y(t),\nabla V(q_y(t))\big)$. By
\eqref{eq:eq3.26} it is
\[
\mu(t)-\mu(0)=\int_0^t\mu'(\tau)\,\text
d\tau=\int_0^t\left[g\big(-\nabla V(q_y),\nabla V(q_y)\big)+
\mathrm H^V(q_y)[\dot q_y,\dot q_y]\right]\,\text d\tau
\]
then, by Proposition \ref{prop:est1} and nondegeneracy of the
maximum point $x_0$, we see that there exists $\nu>0$ such that
\[
\mu(t)-\mu(0)\le-\nu\int_0^t\frac 12\left|\dot q_y\right|^2\,\text
d\tau=-\nu\int_0^t\big(E-V(q_y(\tau))\big)\,\text d\tau,
\]
and since $\lim\limits_{t\to+\infty}\mu(t)=0$ we have
\[
g\big(\dot q_y(0),\nabla V(y)\big)=\mu(0)\ge \nu\int_0^{+\infty}
\big(E-V(q_y(\tau))\big)\,\text d\tau.
\]
Now, consider the map $\kappa(t)=E-V(q_y(t))$: it is
\[
\kappa''(t)=-\mathrm H^V(q_y)[\dot q_y,\dot q_y]+g\big(\nabla
V(q_y),\nabla V(q_y)\big).
\]
Again, by nondegeneracy of $x_0$ as maximum point and Proposition
\ref{prop:est1} there exists $A>0$ such that
\[
g\big(\nabla V(q_y(t)),\nabla V(q_y(t))\big)\le A(E-V(q_y(t)))
\]
while the conservation law of the energy for $q_y$ gives $\frac 12
g\big(\dot q_y,\dot q_y\big)=E-V(q_y)$. Then there exists $B>0$
such that $\kappa''(t)\le B\kappa(t)$ for $t\ge 0$, and by Remark
\ref{rem:rev}
\[
E-V(q_y(t))\ge \big( E- V(y)\big) e^{-\sqrt B t}.
\]
Then
\[
g\big(\dot q_y(0),\nabla V(y)\big)\ge
\nu\big(E-V(y)\big)\int_0^{+\infty}e^{-\sqrt B \tau}\,\text d\tau.
\]
Finally, by Proposition \ref{prop:prop3.14}, $\text d\dot
q_y(0)\to \text d\dot q_{x_0}(0)$ while $\text d\dot
q_{x_0}(0)[v]=\dot\xi_0(0)$ where $\xi_0(t)$ is the unique
solution of
\[
\left\{
\begin{aligned}
&\ddot\xi_0+\mathrm L^V(x_0)[\xi_0]=0\\
&\xi_0(0)=v\\
&\lim\limits_{t\to+\infty}\xi_0(t)=\lim\limits_{t\to+\infty}\dot\xi_0(t)=0.
\end{aligned}
\right.
\]
But, denoting by $\mathbf{e_i}$ a basis of eigenvectors for
$\mathrm L^V(x_0)$ and by $\lambda_i<0$ the corresponding
eigenvalues we have
\[
\xi_0(t)=\sum_{i=1}^m v_i e^{-\sqrt{\lambda_i} t}\mathbf{e_i}.
\]
Since $\text d\dot q_{x_0}(0)[v]=\dot\xi_0(0)$ and $-\mathrm
H^V(x_0)$ is positive definite, there exists $\mu_0>0$ such that
\[
g\big(\text d\dot q_0(0)[v],v\big)\ge\mu_0\,g\big(v,v\big),
\]
and \eqref{eq:eq3.33} is completely proved.
\end{proof}

Finally, we give the result needed to prove our multiplicity
result for homoclinics in \cite{GiaGiaPicI}.
To this aim, take $y\in\{x\,:\,V(x)<E\}$ and consider
\begin{equation}\label{eq:dj}
d(y)=\dist_E(y,V^{-1}(E))
\end{equation}
where $\dist_E$ is the distance with respect to the Jacobi metric.
Combining the results of Theorem~\ref{thm:lem8.5}, Lemma
\ref{thm:lem9.2}, Propositions
\ref{prop:prop3.15}--\ref{prop:prop3.16} and using the function
\eqref{eq:dj} gives us the following:

\begin{teo}\label{thm:propfin}
Assume that:
\begin{itemize}
\item[(a)] $V^{-1}\big(\left]-\infty,E\right[\big)\bigcup\{x_0\}$ is homeomorphic
to an open ball of $\R^m$;
\item[(b)] $\mathrm dV(x)\ne0$ for all $x\in V^{-1}(E)\setminus\{x_0\}$;
\end{itemize}
moreover, let $d$ be as in \eqref{eq:dj}.
Then, there exists a positive number $\delta_*$ such
that, setting
\[
\Omega_*=\{x\in\R^M\,:\,d(x)>\delta_*\}
\]
and denoting by $D_0$ the connected component of
$\partial\Omega_*$ close to $x_0$ and by $D_1$ the connected component of
$\partial\Omega_*$ near $V^{-1}(E)\setminus\{x_0\}$, the following results hold:
\smallskip

\begin{enumerate}
\item $\partial\Omega_*$ is of class $C^2$;
\smallskip

\item
$\overline{\Omega_*}$ is homomorphic to an annulus;
\smallskip

\item
$\overline{\Omega_*}$ is strongly concave with respect to the
Jacobi metric $g_E$;
\smallskip

\item if $x:[0,1]\to\overline{\Omega_*}$ is
an orthogonal geodesic chord in $\overline{\Omega_*}$ relatively
to the Jacobi metric $g_E$ such that $x(0)\in D_0$ and $x(1)\in D_1$, then there
exists $\left]\alpha,\beta\right[\supset [0,1]$ and a unique extension
$\widehat x:[\alpha,\beta]\to\overline\Omega$, $x\in C^0\cap
H^1_{\text{loc}}\big([\alpha,\beta],\overline{\Omega}_E\big)$ satisfying
\smallskip

\begin{itemize}
\item $\widehat x$ is a geodesic with respect to the Jacobi
metric;
\smallskip

\item $\widehat x(s)\in d^{-1}\big(\left]-\delta_*,0\right[\big)$ for all
$s\in \left]\alpha,0\right[\bigcup \left]1,\beta\right[$;
\smallskip

 \item $\widehat x(\alpha)=x_0$,
$\widehat x(\beta)\in V^{-1}(E)\setminus\{x_0\}$.
\end{itemize}
\end{enumerate}
\end{teo}

\end{section}


\begin{thebibliography}{99}

\bibitem{AmbCot} A. Ambrosetti, V. Coti Zelati, {\em Multiple Homoclinic
Orbits for a Class of Conservative Systems}, Rend.\ Sem.\ Mat.\
Univ.\  Padova, Vol.\ 89 (1993), 177--194.

\bibitem{benci} V.\ Benci, {\em Closed Geodesics for the Jacobi Metric
and Periodic Solutions of Prescribed Energy of a Natural
Hamiltonian System}, Ann.\ Inst.\ H.\ Poincar\'e -- Analyse non
Lin\'eaire {\bf1} (1984), 401--412.


\bibitem{bos} W.\ Bos, {\em Kritische Sehenen auf Riemannischen
Elementarraumst\"ucken}, Math.\ Ann.\ {\bf 151} (1963), 431--451.


\bibitem{docarmo} M.\ P.\ do Carmo, {\em Riemannian Geometry}, Birkh\"auser,
Boston, 1992.

\bibitem{GiaGiaPicI} R. Giamb\`o, F. Giannoni, P. Piccione,
{\em Multiple brake orbits and homoclinics in Riemannian manifolds}, in preparation.

\bibitem{GiaGiaPicII} R. Giamb\`o, F. Giannoni, P. Piccione,
{\em Multiple brake orbits   and the Seifert's Conjecture}, in preparation.

\bibitem{gluckziller} H.\ Gluck, W.\ Ziller, {\em Existence
of Periodic Motions of Conservative Systems}, in ``Seminar on
Minimal Surfaces'' (E.\ Bombieri Ed.), Princeton University Press,
65--98, 1983.

\bibitem{Long} Y. Long, Yiming, C. Zhu, {\em Closed characteristics on
compact convex hypersurfaces in $\R\sp{2n}$} Ann.\ of Math.\ (2)
{\bf155}  (2002),  no.\ 2, 317--368.

\bibitem{LustSchn} L.\ Lusternik, L.\ Schnirelman, {\em Methodes
Topologiques dans les Problemes Variationelles}, Hermann, 1934.

\bibitem{Rab} P. H. Rabinowitz, {\em Periodic and Eteroclinic Orbits
for a Periodic Hamiltonian System}, Ann.\ Inst.\ H.\ Poincar\'e,
Analyse Non Lineaire {\bf6} (1989), 331--346.

\bibitem{seifert} H.\ Seifert, {\em Periodische Bewegungen Machanischer
Systeme}, Math.\ Z.\ {\bf51} (1948), 197--216. 

\bibitem{Tan} K. Tanaka, {\em A Note on the Existence of Multiple Homoclinic
Orbits for a Perturbed Radial Potential}, No.\ D.\ E.\ A.\ {\bf1} (1994),
149--162.

\bibitem{Wein} A. Weinstein, {\em Periodic  orbits for convex Hamiltonian
systems}, Ann.\ of Math.\ {\bf108} (1978), 507--518.

\end{thebibliography}
\end{document}